\documentclass[11pt]{article}

\usepackage{amsmath}
\usepackage{amssymb}
\usepackage{amsfonts}
\usepackage{amsthm}
\usepackage{bm}
\usepackage{graphicx}
\usepackage{subcaption}
\usepackage{multirow}
\usepackage{authblk}
\usepackage{natbib}
\usepackage{algorithm}
\usepackage{algorithmic}
\usepackage{dsfont}

\usepackage{tikz}
\usetikzlibrary{arrows,automata}
\usetikzlibrary{shapes.geometric, arrows}
\usetikzlibrary{positioning}

\usepackage{hyperref}
\usepackage[latin1]{inputenc}

\newtheorem{definition}{Definition}
\newtheorem{proposition}{Proposition}

\begin{document}
\title{Improved $x$-space Algorithm for Min-Max Bilevel Problems with an Application to Misinformation Spread in Social Networks}

\author[1]{K\"{u}bra Tan{\i}nm{\i}\c{s}}
\author[2]{Necati Aras}
\author[2]{\.{I}.Kuban Alt{\i}nel}

\affil[1] {Institute of Production and Logistics Management, Johannes Kepler University, Linz, Austria}
\affil[1] {kuebra.taninmis\_ersues@jku.at}
\affil[2] {Department of Industrial Engineering, Bo\u{g}azi\c{c}i University, \.{I}stanbul, Turkey}
\affil[2] { \{arasn,altinel\}@boun.edu.tr}

\date{}

\maketitle

\begin{abstract}
In this work we propose an improvement of the $x$-space algorithm developed for solving a class of min--max bilevel optimization problems (Tang Y., Richard J.P.P., Smith J.C. (2016), A class of algorithms for mixed-integer bilevel min--max optimization. Journal of Global Optimization, 66(2), 225--262). In this setting, the leader of the upper level problem aims at restricting the follower's decisions by minimizing an objective function, which the follower intends to maximize in the lower level problem by making decisions still available to her. The $x$-space algorithm solves upper and lower bound problems consecutively until convergence, and requires the dualization of an approximation of the follower's problem in formulating the lower bound problem. We first reformulate the lower bound problem using the properties of an optimal solution to the original formulation, which makes the dualization step unnecessary. The reformulation makes possible the integration of a greedy covering heuristic into the solution scheme, which results in a considerable increase in the efficiency. The new algorithm referred to as the improved $x$-space algorithm is implemented and applied to a recent min--max bilevel optimization problem that arises in the context of reducing the misinformation spread in social networks. It is also assessed on the benchmark instances of two other bilevel problems: zero-one knapsack problem with interdiction and maximum clique problem with interdiction. Numerical results indicate that the performance of the new algorithm is superior to that of the original algorithm, and also compares favorably with a recent algorithm developed for mixed-integer bilevel linear programs.


\end{abstract}

\section{Introduction and Related Works} \label{section:Introduction}
We address a two-player sequential game, where the second player (follower) has full knowledge of the first player's (leader's) action, and the leader anticipates the follower's optimal reaction in his/her decision. In other words, the leader and follower play a Stackelberg game \citep{stackelberg1952theory}. The leader aims to minimize a function that the follower aims to maximize, i.e., hampers the follower's objective over a feasible solution set that includes the optimal solutions of the follower's maximization problem. This problem can be formulated as the following bilevel optimization model:
\begin{align}
& z^* = \min_{\mathbf{w} \in\mathcal{W}} \, \max_{\mathbf{x}\in \mathcal{X}(\mathbf{w})} \, z(\mathbf{x}).
\end{align}
Here $\mathbf{w}$ and $\mathbf{x}$ represent the decision variables controlled by the leader and follower, respectively, while $\mathcal{W}$ and $\mathcal{X}(\mathbf{w})$ denote their corresponding feasible regions. We consider in particular the case where both the leader and follower make binary decisions.

Bilevel programming problems are $\mathcal{NP}$-hard even for case when upper and lower level problems are linear programs \citep{jeroslow1985polynomial,bard1991some}. When some of the decision variables are restricted to be integers, the problem becomes a mixed-integer bilevel linear program (MIBLP). If integer variables exist only in the upper level problem, then a typical solution approach is to reformulate the problem as a single-level non-linear program using the optimality conditions on the follower's problem. On the other hand, if some of the follower's variables are integer, the standard use of optimality conditions is not possible and solving the bilevel problem becomes generally more difficult. One of the first methods for general MIBLPs is the implicit enumeration scheme due to \cite{moore1990mixed}. \citet{bard1992algorithm} also propose a branch-and-bound algorithm for the case where all variables are binary and use this algorithm to solve problems with up to 50 variables. \citet{DenegreRalphs09} address problems with only integer variables and devise a branch-and-cut method which generates cuts that eliminate bilevel infeasible solutions. \citet{xu2014exact} develop another branch-and-bound method and provide computational results on more than 100 test instances with different sizes up to 920 variables. The cutting plane and branch-and-cut algorithms of \citet{caramia2015enhanced} also eliminate bilevel infeasible solutions as done by \citet{DenegreRalphs09}. They address pure integer  problems and use different valid inequalities which they show to be more effective in terms of solution time on test instances with up to 25 variables. \citet{Fischettietal17} and \cite{wang2017watermelon} use intersection cuts and multi-way disjunction cuts, respectively, to eliminate bilevel infeasible solutions. The branch-and-cut algorithm in the former is considered as the state-of-the-art approach. \citet{lozano2017value} propose a sampling algorithm for general MIBLPs based on the solution of approximate integer upper and a mixed-integer lower level subproblems that are generated using samples of feasible solutions.

A particular class of studies which deal with interdiction problems formulated as MIBLPs where the aim of one player is to damage the objective value of the other's by restricting its decisions, and develops algorithms by exploiting their special structure. A significant amount of these studies deals with the interdiction of some network components by an agent who aims to destroy the functionality of the system. \citet{wood1993deterministic} develops integer programming formulations for a maximum flow interdiction problem, in the context of reducing drug flow in South America. \cite{church2004identifying} introduce the $r$-Interdiction Median Problem (RIM) and $r$-Interdiction Covering Problem (RIC) in which $r$ existing facilities are interdicted with the objective being the maximization of the weighted distances between demand points and their assigned facilities in the RIM and the maximization of the coverage reduction in the RIC. They formulate the problems as mixed-integer linear programs (MILPs) by using closest assignment constraints. \cite{scaparra2008bilevel,scaparra2008exact} add a fortification level to the RIM to obtain the so-called RIM with fortification (RIMF) where the aim is to identify the best protection strategy to reduce the impact of the most disruptive attack. In \cite{scaparra2008exact}, the RIMF is formulated as a single-level maximal covering problem with precedence constraints, whereas in \cite{scaparra2008bilevel} an MIBLP is proposed with the RIM being the lower-level problem. Extended versions of RIMF have been addressed in studies such as \cite{liberatore2011analysis} and \cite{aksen2010budget}.

\cite{Hemmati} address an interdiction problem within the context of influence maximization, and present a cutting-plane algorithm that can optimally solve instances with up to 21 nodes. \cite{tang2016class} develop exact algorithms for general interdiction problems with a mixed-integer follower's problem. Their algorithms keep generating subproblems using samples of feasible solutions, and solving them to compute lower and upper bounds on the optimal value until they converge. The sampling approach is also used by \cite{lozano2016backward} for the three-level defender-attacker-defender games, which they call interdiction problems with fortification. Solutions are sampled from the third level's feasible region and an iterative algorithm is proposed. \cite{fischetti2019interdiction} develop a branch-and-cut algorithm for interdiction problems whose lower level problem satisfies the down-monotonicity assumption which is also necessary for one of the general MIBLP algorithms proposed in \citet{tang2016class}. This assumption restricts the applicability of the algorithm since it is violated in some important applications including the Shortest Path Interdiction Problem \citep{israeli2002shortest} and the Misinformation Spread Minimization Problem \citep{kubra2019misinformation}. The $x$-space algorithm devised in \cite{tang2016class}, however, which is shown to outperform the other two algorithms in the same paper, does not make this assumption, and becomes the starting point for the research in the present paper.

This paper's contribution is twofold. The first one is methodological, and makes algorithmic improvements for the $x$-space algorithm of \citet{tang2016class}. This algorithm is developed for the solution of mixed-integer bilevel min--max optimization problems when all the decision variables of the leader and the integer variables of the follower are restricted to be binary. It solves approximating subproblems to generate lower and upper bounds at each step. One of its key features is the sampling approach used during the creation of the approximations: sample vectors are collected only from the follower's feasible solution set, i.e., $x$-space, resulting in a linear programming (LP) approximation of the second level problem. This LP is then dualized to obtain a single-level approximation whose solution generates a lower bound on the optimal value of the bilevel problem. As pointed out by the authors, most of the computational effort is spent on the generation and solution of the lower bound problems. Hence, any improvement in the efficiency of this step can considerably decrease the solution time of the overall algorithm. We have developed an alternative formulation of the lower bound problem so that the dualization of the follower's problem is no longer necessary. This has been possible by means of \textit{blocking} the follower's solutions. A similar approach is also partly considered for the fortification-interdiction-recourse problems in \cite{lozano2016backward}. However, while they aim to make some interdiction decisions infeasible in the fortification step, we try to determine whether each follower decision (i.e., recourse) is infeasible or not, and ensure that it is blocked in the interdiction step, if necessary. In addition, they do not propose a general procedure for the solution of the restricted interdiction problem in their algorithm. The idea of blocking the follower's feasible solutions is also considered in \cite{lozano2017value}, where a single-level relaxed formulation with blocking constraints on the follower's objective value is solved at every step and new follower solutions are utilized to add constraints to this formulation. The main difference with the improved $x$-space algorithm considered in this paper is that the special structure of the interdiction problems allows a set covering formulation as opposed to a more general and complex mathematical model.

Note that the relation between interdiction and covering problems has already been addressed in the literature. \citet{dinitz2013packing} use linear programming duality to show the relation between the packing interdiction problem with a continuous lower level and the partial covering problem in the sense that the same approximation algorithm can be utilized for both problems. \citet{scaparra2008exact} formulate a fortification-interdiction problem called RIMF as a maximal covering problem by enumerating the interdiction patterns. To the best of our knowledge, this study is the first one proposing an exact algorithm for a general class of discrete bilevel interdiction problems by means of a min-max covering formulation. This feature also helps the improved $x$-space algorithm to benefit from the computational efficiency of covering heuristics.

The second contribution is the application of the improved $x$-space algorithm to exactly solve the Misinformation Spread Minimization Problem (MSMP) for which only heuristics are proposed so far \citep{kubra2019misinformation}. In this problem, which is defined in the context of a social network, the leader protects $h$ nodes, and the follower subsequently activates $k$ nodes among unprotected ones to start a diffusion process. The objective of the leader is to minimize the expected final spread, whereas that of the follower is to maximize the same quantity. The problem has important applications such as preventing the diffusion of misinformation or fake-news as well as determining the vaccination/immunization strategies to prevent the spread of a disease. In fact, it is a stochastic interdiction problem where the leader interdicts (prevents) the activation of a set of nodes by the follower. The bilevel formulation proposed in \cite{kubra2019misinformation} does not have the required structure for the design of the improved algorithm. Nevertheless, the application becomes possible via reformulation and utilizing some problem specific features. Note that the MSMP does not satisfy the down-monotonicity assumption, and thus the general MIBLP solution algorithm by \cite{fischetti2019interdiction} is not applicable.

The remainder of the paper is organized as follows. Section \ref{section:x-space} outlines the original $x$-space algorithm, and discusses the properties of the optimal solutions to the lower bound problem. The modifications to improve the original algorithm are also described and analyzed in this section. In Section \ref{section:MSMP}, we define the MSMP and modify the current formulation with the aim of implementing the improved  $x$-space algorithm. Section \ref{section:Numerical Results} is devoted to the numerical experiments. We conclude the paper in Section \ref{section:Conclusion} by mentioning some potential future research directions.

\section{Improved $x$-space Algorithm} \label{section:x-space}

\subsection{The Original Algorithm} \label{subsection:original_algorithm}
\cite{tang2016class} focus on problems that can be addressed in the form of a Stackelberg game where both the leader and the follower have the same objective function. The leader tries to minimize it, whereas the follower's aim is its maximization. They define the bilevel optimization problem under investigation as
\begin{equation}
Z^\ast(\mathcal{W},\mathcal{F})=\min_{\mathbf{w}\in\mathcal{W}}\max_{\mathbf{x}\in\mathcal{X}}\left\{\mathbf{p}^T \mathbf{x} : (\mathbf{w},\mathbf{x})\in \mathcal{F}\right\}, \label{min_max problem}
\end{equation}
where $\mathbf{w}\in \mathbb{Z}^{n_1}$ denotes the upper level decision variables and $\mathbf{x} = (\mathbf{x}_1,  \mathbf{x}_2)^T$ with $\mathbf{x}_1 \in \mathbb{R}^{n_2-q}$ and $\mathbf{x}_2\in\mathbb{Z}^q$ denoting the lower level decision variables. As can be seen, the leader's decision variables are discrete while the follower's decision variables can be both continuous and discrete. The feasible solution set $\mathcal{F}$ is defined as
\begin{equation}
\mathcal{F}=\left\{\mathbf{w}\in \mathcal{W}, \mathbf{x}\in \mathcal{X}: \mathbf{C}\mathbf{x} + \mathbf{D} \mathbf{w}\leq \mathbf{d} \right\},
\end{equation}
where $\mathcal{W}=\{\mathbf{w}\in \mathbb{Z}^{n_1}:\mathbf{A_1} \mathbf{w}\leq \mathbf{b_1}, \mathbf{0}\leq \mathbf{w}\leq \mathbf{u}\}$ denotes the feasible region of the upper level problem, and the feasible region of the lower level problem is defined using the following two sets of constraints: $\mathcal{X}=\{\mathbf{x}_1 \in \mathbb{R}^{n_2-q},\mathbf{x}_2\in\mathbb{Z}^q: \mathbf{A_2} \mathbf{x}\leq \mathbf{b_2}, \mathbf{x}\geq \mathbf{0}\}$ is independent of the upper level decisions and $\mathbf{C}\mathbf{x} + \mathbf{D} \mathbf{w}\leq \mathbf{d}$ relates the upper and lower level variables. Here, $\mathbf{A_1}$, $\mathbf{A_2}$, $\mathbf{C}$, and $\mathbf{D}$ are $q_1 \times n_1$, $q_2 \times n_2$, $q_3 \times n_2$, and $q_3 \times n_1$ dimensional rational matrices, respectively, while $\mathbf{u}$, $\mathbf{b_1}$, $\mathbf{b_2}$, and $\mathbf{d}$ are vectors of size $n_1$, $q_1$, $q_2$, and $q_3$, respectively.

Let $\mathcal{W}^\circ=\mathcal{W}\cap\{0,1\}^{n_1}$, $\mathcal{X}^\circ=\mathcal{X}\cap\{0,1\}^{n_2}$, and $\mathcal{F}^\circ=\mathcal{F}\cap\{0,1\}^{n_1+n_2}$. For $\mathbf{C}=\mathbf{I}$ and $\mathbf{u}=\mathbf{1}$, it is shown that
\begin{equation}
Z^{\ast}(\mathcal{W}^\circ,\mathcal{F}^\circ)=Z^\ast(\mathcal{W}^\circ,\widetilde{\mathcal{F}^\circ}),
\end{equation}
where
$ \widetilde{\mathcal{F}^\circ}=\left\{ \mathbf{w}\in \mathbb{R}^{n_1}, \mathbf{x}\in \mathbb{R}^{n_2} : \mathbf{w}\in \mathcal{W}^\circ, \mathbf{x}\in \text{conv($\mathcal{X}^\circ$)}, \mathbf{x}\leq \mathbf{d}- \mathbf{D} \mathbf{w} \right \} $ \cite[p.~242]{tang2016class}. Notice that the follower's problem is an LP when the feasible region is defined by $\widetilde{\mathcal{F}^\circ}$, in which case it is possible to reformulate \eqref{min_max problem} as a single-level problem. Nevertheless, to obtain the convex hull of $\mathcal{X}^\circ$ is usually intractable. Therefore, \cite{tang2016class} propose an algorithm that involves sampling from the solution space of the lower level problem (which they call the $x$-space) and obtaining an (inner) approximation of the conv($\mathcal{X}^\circ$) iteratively. They also show that this approximation leads to a lower bound on the optimal objective value $Z^{\ast}(\mathcal{W}^\circ,\mathcal{F}^\circ)$.

The original $x$-space algorithm, which is given as Algorithm \ref{algorithm:XS} below, requires the following assumptions:
\begin{enumerate}
\item The decision variables of the leader's problem are binary.
\item The feasible region $\mathcal{W}^\circ$ of the leader's problem is not empty.
\item The follower's problem is feasible for each decision $\mathbf{w}\in \mathcal{W}^\circ$ of the leader.
\item The follower's problem is bounded above for at least one decision of the leader.
\item $\mathbf{x}=\mathbf{0}$ is a feasible follower's solution for each leader decision $\mathbf{w} \in \mathcal{W}^\circ$.
\end{enumerate}

\begin{algorithm}[h]
\caption{The $x$-space algorithm} \label{algorithm:XS}
\begin{algorithmic}[1]
\STATE {\textbf{Step 0: Initialization} }
\STATE {\indent Define $J=\{-(\rho-1),...,0\}$ as an index set of initial solutions where $\rho \in \mathbb{Z}_+. $}
\STATE {\indent Select $\mathbf{x}^0$ such that $\mathbf{x}^0$ is a feasible follower solution for all feasible leader solutions.}
\STATE {\indent Choose a feasible $\mathbf{w}^\tau$ and  $\mathbf{x}^\tau \in \arg\max_\mathbf{x} \{ \mathbf{p}^T \mathbf{x} : (\mathbf{w}^\tau,\mathbf{x})\in \mathcal{F}^\circ \big\}$ for $\tau \in J \setminus \{0\}$.}
\STATE {\indent Let $S_1=\bigcup_{\tau \in J}\{\mathbf{x}^\tau\}$.}
\STATE {\indent Set $q=1$.}
\STATE {\textbf{Step $\mathbf{q}$: }$(q=1,2,...)$}
\STATE {\hspace{10mm}\textbf{Step $\mathbf{q_1}$:} Obtain an optimal solution $\mathbf{w}^q$ and optimal value $l_q$ to the approximate leader problem $(LB_q)$.}
\STATE {\hspace{10mm} \textbf{Step $\mathbf{q_2}$:} Obtain an optimal solution $\mathbf{x}^q$ and optimal value $u_q$ to the follower problem $(UB_q)$.}
\STATE {\hspace{10mm} \textbf{Step $\mathbf{q_3}$:} \textbf{if} $l_q\neq u_q$ \textbf{then}}
\STATE {\hspace{35mm} Expand $S_{q+1}=S_q\cup \{\mathbf{x}^q\}$ and update $J=J\cup {q}$.}
\STATE {\hspace{35mm} \textbf{Go to} Step $q+1$.}
\STATE {\hspace{28mm} \textbf{else}}
\STATE {\hspace{35mm} Return optimal solution $(\mathbf{w}^*,\mathbf{x}^*)=(\mathbf{w}^q,\mathbf{x}^q)$ and $z^*=u_q$.}
\end{algorithmic}
\end{algorithm}

The upper bound problem $UB_q$ in Algorithm \ref{algorithm:XS} is simply the follower's problem which is defined using the most recently obtained leader solution $\mathbf{w}^q$:
\begin{equation}
u_q=\max_\mathbf{x}{}\{\mathbf{p}^T \mathbf{x} : (\mathbf{w}^q,\mathbf{x})\in \mathcal{F}^\circ\} \geq Z^{\ast}(\mathcal{W}^\circ,\mathcal{F}^\circ).
\end{equation}

The lower bound problem $LB_q$ given as
\begin{align}
l_q =&\min_{\mathbf{w}\in \mathcal{W}^\circ}\max_{\mathbf{x}} \{\mathbf{p}^T \mathbf{x} : \mathbf{x} \in \text{conv}(S_q), \mathbf{x} \leq \mathbf{d}-\mathbf{D} \mathbf{w}\} \notag \\
=& \min_{\mathbf{w}\in \mathcal{W}^\circ}\max_{\mbox{\scriptsize\boldmath$\lambda$},\mathbf{x}}
\left\{ \mathbf{p}^T \mathbf{x} : \sum_{\tau \in J} \lambda_\tau=1, \mathbf{x}-\sum_{\tau\in J}\lambda_\tau \mathbf{x}^\tau=\mathbf{0}, \mathbf{x} \leq \mathbf{d}- \mathbf{D} \mathbf{w}, \mathbf{x}\geq \mathbf{0},\, \lambda_\tau\geq 0,\,  \tau\in J  \right\}
 \label{equation:LB problem}
\end{align}
is a bilevel problem itself, except that $\mathbf{x}\in \text{conv($\mathcal{X}^\circ$)}$ is replaced by $\mathbf{x}\in \text{conv($S_q$)}$, where $S_q$ approximates the follower's solution set at step $q$. The set $J$ contains the indices of the solutions in $S_q$, and the new decision variables $\lambda_{\tau}$ are defined in order to express the follower solution $\mathbf{x}$ as a convex combination of the solutions $\mathbf{x}^\tau$ in $S_q$ using the constraint $\mathbf{x}-\sum_{\tau\in J}\lambda_\tau \mathbf{x}^\tau=\mathbf{0}$. Together with $\sum_{\tau \in J} \lambda_\tau=1$, this constraint implies that $\mathbf{x}\in \text{conv}(S_q)$. Notice that the constraint $\mathbf{x} \leq \mathbf{d}- \mathbf{D} \mathbf{w}$, which links the upper and lower level variables, is included in $LB_q$ formulation as opposed to the set $\mathcal{X}$ whose convex hull is inner-approximated. \cite{tang2016class} show that the lower bound obtained at each iteration is at least as good as the previous one and the algorithm converges to an optimal solution in finite number of iterations.

Notice that $LB_q$'s lower level problem is an LP. It is solved by taking the dual of the LP by treating the leader's variable vector $\mathbf{w}$ as a parameter in \cite{tang2016class}. Then, the constraints describing $\mathcal{W}^\circ$ are added to this formulation to end up with a non-linear model since $\mathbf{w}$ is actually a decision variable. After removing the nonlinearities with the help of newly defined decision variables necessary for linearization, the final single-level MILP is solved by a commercial solver. The correctness of this scheme is based on the fact that the leader and follower objective functions are the same but optimized in opposite directions. In the following sections, we propose a more efficient method to optimally solve $LB_q$. As we have mentioned in the introduction, \cite{tang2016class} point out that more efficient solution of $LB_q$ is one of the future research challenges that can markedly improve the performance of the $x$-space algorithm.

\subsection{A New Formulation for the Lower Bound Problem} \label{subsection:new_LB}

We focus our attention on the setting with  $\mathbf{d}=\mathbf{1}$, $\mathbf{C}=\mathbf{I}$, $\mathbf{D}=\mathbf{I}$, and $n_1=n_2=n$ in the development of the alternative formulation. It will be clear in Section \ref{section:Implementation} and Section \ref{subsection:Comparison} that the method can be adapted to problems with other parameter values for $\mathbf{d}$, $\mathbf{C}$, and $\mathbf{D}$, and different number of decision variables in the upper and lower level problems as long as some or all of the follower's decisions are interdicted by the leader.

Let $N=\{1,...,n\}$ denote the set of indices for both upper and lower level variables. First, we consider the inner optimization problem of $LB_q$ in \eqref{equation:LB problem}, which we denote by $LB_q(\mathbf{w})$. Notice that $\mathbf{w}$ is not a decision variable but a parameter for $LB_q(\mathbf{w})$. Given $\mathbf{w} \in \mathcal{W}^\circ$ and $S_q$, we define the following set
\begin{equation}
B(\mathbf{w},S_q)=\{\tau \in J : \exists i\in N, x_i^\tau > 1-w_i \},
\end{equation}
and its complement $\Bar{B}(\mathbf{w},S_q)=J \setminus B(\mathbf{w},S_q)$. By its definition, $B(\mathbf{w},S_q)$ contains the indices of the follower's solutions in $S_q$ that are rendered infeasible, i.e., \textit{blocked}, by the leader's solution $\mathbf{w}$. The remaining ones are not blocked and are represented by $\Bar{B}(\mathbf{w},S_q)$.

\begin{proposition}
Let $z_\tau$ denote the objective value of the follower's solution $\mathbf{x}^\tau$. The optimal objective value of $LB_q(\mathbf{w})$ is
\begin{equation}
z^\ast(\mathbf{w},S_q)=\max_{\tau \in \bar{B}(\mathbf{w},S_q)} z_\tau.
\end{equation}
\label{Prop:optimal_obj}
\end{proposition}

\begin{proof}
Firstly, we will show that if $(\mbox{\boldmath$\lambda$},\mathbf{x})$ is a feasible solution to $LB_q(\mathbf{w})$, then $\lambda_\tau=0$ for all $\tau \in B(\mathbf{w},S_q)$.
Since $(\mbox{\boldmath$\lambda$},\mathbf{x})$ is feasible, $x_i=\sum_{\tau\in J}\lambda_\tau x_i^\tau\leq 1- w_i$ must hold for all $i\in N$. By contradiction, suppose that there exists a $\tau^\prime \in B(\mathbf{w},S_q)$ such that $\lambda_{\tau^\prime}>0$. Then, there must exist at least one $i\in N$ such that $w_i=1$ and $x_i=\sum_{\tau\in J}\lambda_\tau x_i^\tau>0$, which contradicts $x_i \leq 1-w_i$.

Note that the objective function of the second level maximization problem can be arranged as
\begin{equation}
\mathbf{p}^T \mathbf{x}=\sum_{i\in N} p_i x_i=\sum_{i\in N}p_i \sum_{\tau\in J}\lambda_\tau x_i^\tau=\sum_{\tau\in J}\lambda_\tau \sum_{i\in N}p_ix_i^\tau.
\end{equation}
The objective value of the follower's solution $\mathbf{x}^\tau$ is defined as $z_\tau=\sum_{i\in N}p_ix_i^\tau$. Thus, $\mathbf{p}^T \mathbf{x}=\sum_{\tau\in J}\lambda_\tau z_\tau$ follows.
Using the result that $\lambda_\tau=0$ for all $\tau \in B(\mathbf{w},S_q)$, the objective function of $LB_q(\mathbf{w})$ can be rewritten as $\mathbf{p}^T \mathbf{x}=\sum_{\tau\in J}\lambda_\tau z_\tau=\sum_{\tau\in \bar{B}(\mathbf{w},S_q)}\lambda_\tau z_\tau$.
Similarly, $\mathbf{x}=\sum_{\tau \in J}\lambda_\tau \mathbf{x}^\tau =\sum_{\tau \in \bar{B}(\mathbf{w},S_q)}\lambda_\tau \mathbf{x}^\tau \leq \mathbf{1}-\mathbf{w}$ because we know that $\mathbf{x}^\tau\leq \mathbf{1}-\mathbf{w}$ for all $\tau \in \bar{B}(\mathbf{w},S_q)$ by definition, and $\mbox{\boldmath$\lambda$}\geq \mathbf{0}$. Then, by fixing $\lambda_\tau=0, \tau \in B(\mathbf{w},S_q)$, the problem reduces to
\begin{equation}
z^\ast(\mathbf{w},S_q)=\max_{\mbox{\scriptsize\boldmath$\lambda$}} \left \{\sum_{\tau\in \bar{B}(\mathbf{w},S_q)}\lambda_\tau z_\tau :  \sum_{\tau\in \bar{B}(\mathbf{w},S_q)}\lambda_\tau=1, \mbox{\boldmath$\lambda$}\geq \mathbf{0} \right \}.
\end{equation}
As the objective value is a convex combination of the previous objective values, its optimal value is equal to $\max_{\tau \in \bar{B}(\mathbf{w},S_q)} z_\tau$.

\end{proof}

As a consequence of Proposition \ref{Prop:optimal_obj}, the bilevel lower bound problem \eqref{equation:LB problem} becomes
\begin{align}
l_q &= \min_{\mathbf{w}\in \mathcal{W}} \max_{\tau \in \bar{B}(\mathbf{w},S_q)} z_\tau \\
& =\min_{\mathbf{w},z}\, z\\
& \hspace{15pt}\text{s.t.} \notag \\
& \hspace{25pt}\mathbf{w}\in \mathcal{W} \\
&  \hspace{25pt} z \geq z_\tau , \tau \in \bar{B}(\mathbf{w},S_q). \label{equation:z_z_tau}
\end{align}

Instead of generating the set $B(\mathbf{w},S_q)$ for all $\mathbf{w}\in \mathcal{W}$, we define a new set $C_\tau=\{i\in N : x_i^\tau=1\}$, the indices of the leader variables that cause the solution $\mathbf{x}^\tau$ to be blocked when at least one of them is positive. We refer to $C_\tau$ as the \textit{blocker set} of $\tau$. The following proposition establishes the relation between $C_\tau$ and $B(\mathbf{w},S_q)$.

\begin{proposition}
If $\sum_{i\in C_\tau} w_i \geq 1$, then $\tau \in B(\mathbf{w},S_q)$, i.e., $\mathbf{x}^\tau$ is blocked by $\mathbf{w}$.
\label{Prop:B_C relation}
\end{proposition}

\begin{proof}
If $\sum_{i\in C_\tau} w_i \geq 1$, then there exists $i^\prime \in C_\tau$ such that $w_{i^\prime}=1$. By definition of $C_\tau$, $\mathbf{x}^\tau_{i^\prime}=1$. Recall that $B(\mathbf{w},S_q)=\{\tau \in J : \exists i\in N, x_i^\tau > 1-w_i \}$. Since $w_{i^\prime}=1$ and $x^\tau_{i^\prime}=1$, $\tau$ belongs to $B(\mathbf{w},S_q)$.
\end{proof}

Notice that Constraint \eqref{equation:z_z_tau} is stated only for the unblocked solutions. We now introduce a new binary decision variable $\alpha_\tau$ for each $\tau \in J$, to indicate whether a follower solution is blocked or not. In the new formulation denoted as $LB_q'$, $\alpha_\tau=0$ when $\mathbf{x}^\tau$ cannot be blocked and Constraint \eqref{equation:obj_t} becomes $z \geq z_\tau $.
\begin{align}
LB_q':
l_q =& \min z \\
& \text{s.t.} \notag\\
& \hspace{10pt}\mathbf{w}\in \mathcal{W} \\
& \hspace{10pt} z\geq (1-\alpha_\tau)z_\tau &\tau& \in J \label{equation:obj_t}\\
& \hspace{10pt} \alpha_\tau \leq \sum_{i\in C_\tau} w_i &\tau& \in J\\
& \hspace{10pt} z\geq 0, \, \alpha_\tau \in \{0,1\} &\tau& \in J.
\label{equation:LB problem-2}
\end{align}

The new formulation can be improved further using the information on the blocked solutions. The following proposition yields this stronger formulation.
\begin{proposition}
Let $\bar{Z}$ be an upper bound on the optimal objective value $Z^\ast(\mathcal{W},\mathcal{F})$ of the leader. Then, any solution $\mathbf{x}^\tau$ with a larger objective value than $\bar{Z}$ must be blocked in an optimal leader solution, i.e., $\alpha_\tau=1$ for all $\tau$ such that $z_\tau>\bar{Z}$, in an optimal solution of $LB_q'$.
\label{Prop:blockedSolns}
\end{proposition}

\begin{proof}
Let $(\mathbf{w}^*, \mbox{\boldmath$\alpha$}^*)$ be an optimal solution of $LB_q'$. Then $l_q \geq (1-\alpha^*_\tau)z_\tau$ for $\tau\in J$. Suppose there exist $\tau^\prime \in J$ such that $\alpha^*_{\tau^\prime}=0$ and $z_{\tau^\prime}>\Bar{Z}$, which implies $l_q \geq z_{\tau^\prime} > \Bar{Z}$. This clearly contradicts $l_q \leq Z^\ast(\mathcal{W},\mathcal{F}) \leq \bar{Z}$.
\end{proof}

Recall that each solution in $S_q$ can be obtained by solving the follower's problem to optimality with the exception of the trivial solution ($\tau=0$) that is feasible for all feasible leader solutions. Therefore, each $z_\tau, \tau \in J\setminus \{0\}$ constitutes an upper bound on $Z^\ast(\mathcal{W}^\circ,\mathcal{F}^\circ)$, i.e., $z_\tau \geq Z^\ast(\mathcal{W}^\circ,\mathcal{F}^\circ),  \tau \in J\setminus \{0\}$. From now on, $\Bar{Z}=\min_{\tau \in J\setminus \{0\}} z_\tau$ denotes the current upper bound on the optimal objective value. Now, we define $J^B=\{\tau \in J : z_\tau>\Bar{Z}\}$ as the index set of the follower solutions that must be blocked according to Proposition \ref{Prop:blockedSolns}. Then, we can derive
\begin{align}
LB_q'':
l_q =& \min z \label{equation:LB3 obj}\\
& \text{s.t.} \notag \\
& \hspace{10pt}\mathbf{w}\in \mathcal{W} \\
& \hspace{10pt} z\geq (1-\alpha_\tau)z_\tau &\tau& \in J\setminus J^B \label{equation:obj_t2}\\
& \hspace{10pt} \alpha_\tau \leq \sum_{i\in C_\tau} w_i &\tau& \in J\setminus J^B\\
& \hspace{10pt} 1 \leq \sum_{i\in C_\tau} w_i &\tau& \in J^B \label{equation:superValid}\\
& \hspace{10pt} z\geq 0, \, \alpha_\tau \in \{0,1\} &\tau& \in J\setminus J^B
 \label{equation:LB problem-3}
\end{align}
as an equivalent and stronger formulation of $LB_q'$. The constraint \eqref{equation:superValid} is similar to the supervalid inequalities used by \cite{israeli2002shortest}. Notice that the objective value of this formulation is equal to the maximum unblocked follower objective value. In other words, this reformulation is possible since the leader and the follower have the same objective function.


Consider $z_0$, the objective value of trivial follower solution $\mathbf{x}^0$ gives. It is associated with a feasible (but not necessarily optimal) follower reaction for any leader solution, therefore $z_0 \leq \max_\mathbf{x} \{ \mathbf{p}^T \mathbf{x} : (\mathbf{w},\mathbf{x})\in \mathcal{F}^\circ\} \leq Z^\ast(\mathcal{W}^\circ,\mathcal{F}^\circ)$ for all $\mathbf{w} \in \mathcal{W}^\circ$. Furthermore, since $\mathbf{x}^0$ is feasible for all $\mathbf{w} \in \mathcal{W}$, it cannot be blocked by the leader, i.e., $\alpha_0=0$ and $z\geq z_0$ in an optimal solution of $LB_q''$, and $z_0 \leq l_q \leq \bar{Z}$ follows. Next proposition reduces the domain of $l_q$ further and it is a direct consequence of Proposition \ref{Prop:blockedSolns}.

\begin{proposition}
Either $l_q=z_0$ or $l_q=\bar{Z}$ for any $q$.
\label{Prop:two_obj_values}
\end{proposition}

\begin{proof}
First of all, $z_0\leq \bar{Z}$ as explained above. Now, notice that the set $\{z_\tau : \tau \in J\setminus J^B\}$ has only two distinct values, $z_0$ and $\bar{Z}$ since $J\setminus J^B=\{\tau \in J : z_\tau \leq \bar{Z}\}$ and $\Bar{Z}=\min_{\tau \in J\setminus \{0\}} z_\tau$. As a result, Constraint \eqref{equation:obj_t2} includes $z\geq z_0$ and $z\geq (1-\alpha_\tau)\bar{Z}$ for $\tau \in J\setminus (J^B \cup \{0\})$. If there is a feasible leader solution that can block all $\tau \in J\setminus (J^B \cup \{0\})$, then $\alpha_\tau=1$ and $l_q=z_0$; otherwise $\alpha_\tau=0$ and $l_q=\bar{Z}$.
\end{proof}


In other words, if $l_q=\Bar{Z}$, then $\Bar{Z}\leq Z^\ast(\mathcal{W}^\circ,\mathcal{F}^\circ)\leq \Bar{Z}$  and the algorithm terminates; otherwise, the lower bound remains the same as $z_0$.

\subsection{Greedy Maximum Covering}
In an interdiction problem, the upper level feasible region is typically associated with cardinality/budget constraints and logical restrictions on the interdiction decisions. In this part, we focus on the case with cardinality constraints, i.e., $\mathcal{W}^\circ=\{\mathbf{w} \in \{0,1\}^n : \sum_i w_i\leq g\}$. Given the index set $J$ for $S_q$, consider the maximum covering problem
\begin{align}
c_q =& \max \sum_{\tau \in J\setminus \{0\}} \alpha_\tau\\
& \text{s.t.} \notag \\
& \hspace{10pt}\sum_{i\in N} w_i\leq g \\
& \hspace{10pt} \alpha_\tau \leq \sum_{i\in C_\tau} w_i &\tau& \in J\setminus \{0\} \\
& \hspace{10pt}w_i \in \{0,1\},\, \alpha_\tau \in \{0,1\} &i& \in N, \tau \in J\setminus \{0\}
 \label{equation:MC problem}
\end{align}
where the objective is to cover the maximum number among $|J|-1$ items and $\mathbf{x}^{\tau}$ is covered only if $\sum_{i\in C_\tau}w_i\geq 1$. If $c_q= |J|-1$, which is the maximum objective value attainable, then it indicates the existence of a feasible $\mathbf{w}$ which covers (blocks in our context) all $\mathbf{x}^\tau$, $\tau \in J \setminus \{0\}$. Then, the optimum objective value $l_q$ of the lower bound problem for $S_q$ is $z_0$ due to Proposition \ref{Prop:two_obj_values}. Now let $\hat{c}_q$ denote the objective value that a greedy maximum covering algorithm yields for the same problem and $\hat{\mathbf{w}}^q$ denote the corresponding solution. Notice that $\hat{c}_q\leq c_q$ since $\hat{\mathbf{w}}^q$ is a heuristic solution to the problem. If $\Hat{c}_q = |J|-1$, then $\hat{\mathbf{w}}^q$  is an optimal solution to $LB_q$ and $l_q=z_0$. This identity allows us to first identify a good, possibly optimal, feasible  solution of $LB_q$ using a heuristic, and then solve $LB_q''$ given with \eqref{equation:LB3 obj}--\eqref{equation:LB problem-3} if $\hat{c}_q < |J|-1$. At this point, we have all the ingredients of the improved $x$-space algorithm (IXS), which is formally given as Algorithm \ref{algorithm:IXS} in the following. IXS guarantees to find an optimal solution due to the correctness of the original algorithm and the fact that solving $LB_q''$ yields an optimal solution to $LB_q$ as shown before. The usage of the greedy heuristic does not change the correctness of the algorithm since $LB_q''$ is solved optimally as an MILP whenever the greedy algorithm fails to find an optimal solution. It is also worthwhile noticing that the IXS algorithm may yield a different optimal solution from the one provided by the original algorithm in the case that multiple optimal solutions exist for the problem.

\begin{algorithm}[h]
\caption{Improved $x$-space (IXS) algorithm} \label{algorithm:IXS}
\begin{algorithmic}[1]
\STATE {\textbf{Step 0: Initialization} }
\STATE {\indent Define $J=\{-(\rho-1),...,0\}$ as an index set of initial solutions where $\rho \in \mathbb{Z}_+. $}
\STATE {\indent Select $\mathbf{x}^0$ such that it is a feasible follower solution for all feasible leader solutions.}
\FOR {$\tau \in J\setminus\{0\}$ }
\STATE{Choose a feasible $\mathbf{w}^\tau$ and  $\mathbf{x}^\tau \in \arg\max_\mathbf{x} \big\{\mathbf{p}^T \mathbf{x} : (\mathbf{w}^\tau,\mathbf{x})\in \mathcal{F}^\circ \big\}$}
\STATE {Set $z_\tau \leftarrow \mathbf{p}^T \mathbf{x}^\tau$ and $C_\tau \leftarrow \{i\in N : x_i^\tau=1\}$ }
\ENDFOR
\STATE {\indent Let $S_1=\bigcup_{\tau \in J}\{\mathbf{x}^\tau\}$.}
\STATE {\indent Set $q=1$, $\bar{Z}\leftarrow \min_{\tau\in J\setminus\{0\}} z_\tau $ and $\mathbf{w}^*\in \{\mathbf{w}^\tau : z_\tau=\bar{Z}\}$}
\STATE {\textbf{Step $\mathbf{q}$: }$(q=1,2,...)$}
\STATE {\hspace{1cm}\textbf{Step $\mathbf{q_0}$:} Use the Greedy Covering heuristic to obtain $\hat{\mathbf{w}}^q$ and $\hat{c}_q$ }
\STATE {\hspace{2.8cm} \textbf{if} $\hat{c}_q<|J|$ \textbf{then} }
\STATE {\hspace{3.5cm} Go to Step ${q_1}$}
\STATE {\hspace{2.8cm} \textbf{else}}
\STATE {\hspace{3.5cm} Set $\mathbf{w}^q \leftarrow\hat{\mathbf{w}}^q$, $l_q \leftarrow\hat{c}_q$, go to Step $q_2$ }
\STATE {\hspace{1cm}\textbf{Step $\mathbf{q_1}$:} Solve $LB_q''$ to obtain an optimal solution $\mathbf{w}^q$ and optimal value $l_q$}
\STATE {\hspace{1cm}\textbf{Step $\mathbf{q_2}$:} Solve $UB_q$ to obtain an optimal solution $\mathbf{x}^q$ and optimal value $u_q$}
\STATE {\hspace{2.8cm} \textbf{if} $u_q<\bar{Z}$ \textbf{then} }
\STATE{\hspace{3.5cm} Set $\bar{Z}\leftarrow u_q$ and $\mathbf{w}^*\leftarrow \mathbf{w}^q$}
\STATE {\hspace{1cm}\textbf{Step $\mathbf{q_3}$:} \textbf{if} $l_q < \bar{Z}$ \textbf{then}}
\STATE {\hspace{3.5cm} Expand $S_{q+1}=S_q\cup \{\mathbf{x}^q\}$ and update $J=J\cup {q}$.}
\STATE {\hspace{3.5cm} \textbf{Go to} Step $q+1$.}
\STATE {\hspace{2.8cm} \textbf{else}}
\STATE {\hspace{3.5cm} \textbf{Return} optimal solution $(\mathbf{w}^*,\mathbf{x}^*)=(\mathbf{w}^q,\mathbf{x}^q)$ and $z^*=\bar{Z}$.}
\end{algorithmic}
\end{algorithm}
\section{Misinformation Spread Minimization Problem}\label{section:MSMP}

Influence Maximization Problem (IMP) is well-known and defined on social networks with many application areas such as viral marketing. It involves finding a set of $k$ nodes (seed set) to start influence propagation on a network so that the expected number of affected nodes is maximized at the end. \cite{Kempe} define IMP as a stochastic discrete optimization problem for the first time and prove its $\mathcal{NP}$-hardness for widely accepted (stochastic) diffusion models.

\begin{definition}[IMP \citep{Kempe}]
Let $G=(V,A)$ denote a directed graph with node set $V$ and arc set $A$, and $\sigma_M(G,Y)$ denote the expected number of influenced (active) nodes at the end of the diffusion process, i.e. the spread, in graph $G$ under the diffusion model $M$ when the nodes in $Y$ is activated initially. Given $G=(V,A)$ and a positive integer $h$, the IMP is
\begin{equation}
\max_{\substack{Y\subseteq V, \\ |Y|\leq h }} \sigma_M(G,Y).
\label{equation:MSMP}
\end{equation}

\end{definition}

Since then, the IMP and its several versions have been studied extensively \citep{chen2013information}. A recent survey on influence maximization problems can be found in \cite{li2018influence}. \cite{kubra2019misinformation} address the Misinformation Spread Minimization Problem (MSMP), a competitive version of the IMP in the form of a Stackelberg game where the leader of the game protects a subset of nodes and then the follower activates a set of unprotected nodes to start a diffusion process. The aim of the follower is to maximize the expected number of influenced nodes, i.e., the spread, while the leader wants to minimize the same objective function.

\begin{definition}[MSMP \citep{kubra2019misinformation}]
Let $G_{X}$ denote the graph obtained from $G$ by removing the node set $X\subset V$ and all arcs incident to the nodes in $X$. Given a directed graph $G=(V,A)$, a diffusion model $M$, and positive integers $h$ and $k$, the MSMP involves finding a set of nodes $X$ of size $h$ minimizing $\sigma_M(G_X,Y^{*}(X))$, where $Y^{*}(X)$ is a set of $k$ nodes in $V\setminus X$ which maximizes the spread on $G_X$, i.e.,
\begin{equation}
\min_{\substack{X\subset V, \\ |X|\leq h }} \sigma_M(G_X,Y^*(X)),
\label{equation:MSMP}
\end{equation}
where
\begin{equation}
Y^*(X)=\operatorname*{arg\,max}_{\substack{Y\subseteq V\setminus X, \\ |Y|\leq k}}\sigma_M(G_X,Y).
\end{equation}
\end{definition}

\subsection{A Bilevel Programming Formulation via Live-arc Representation}

\cite{kubra2019misinformation} adopt the well-known Linear Threshold (LT) model where each node $i\in V$ has a random threshold value $\theta_i\in(0,1)$, each arc $(i,j)\in A$ has a deterministic weight $w_{ij}$ satisfying $\sum_i w_{ij}\leq 1$, and node $i$ becomes active at any step of the diffusion process if the total incoming weight from its active neighbors exceeds the threshold. Once a node is active, it remains so until the end of the diffusion process. \cite{Kempe} propose the \textit{live-arc} model which is equivalent to the LT under certain assumptions and it leads to a more efficient mathematical formulation. In this model, each diffusion scenario is represented by a directed subgraph of the original directed graph obtained by randomly labeling a subset of the arcs as \textit{live}. If node $i$ is reachable from any seed node via the live arcs in a scenario, then it means that $i$ is influenced in that scenario. For the equivalence of the two models, it is required that the probability of labeling arc $(i,j)$ as live is $w_{ij}$ and for each node at most one of its incoming arcs can be live in each scenario. In the bilevel model formulated using the \textit{live-arc} representation, the follower's problem is a two-stage stochastic program and the uncertainty is handled by realizing a scenario in the form of a directed subgraph of the original graph. The list of sets, parameters and decision variables used in the model are presented below, prior to the formulation.

\begin{tabular}{lllp{12cm}}
\multicolumn{4}{l}{{\textit{Sets and Parameters:}}} \\
&$V$	&:& set of nodes in the network \\
&$h$	&:& number of nodes the leader protects \\
&$k$	&:& number of nodes the follower activates  \\
&$R$	&:& a set of live-arc scenarios \\
&$p_r$		&:& probability of scenario $r$, $r\in R$ \\
&$a_{ji}(\mathbf{w},r)$	&:& 1 if node $i$ is reachable from node $j$ in scenario $r$, when the protection decision is $\mathbf{w}$; 0 otherwise \\
\multicolumn{4}{l}{{\textit{Decision variables:}}}\\
&$w_i$&:& 1 if node $i$ is protected by the leader; 0 otherwise \\
&$y_i$ &:& 1 if node $i$ is activated by the follower; 0 otherwise \\
&$u_{ir}$	&:& 1 if node $i$ is influenced in scenario $r$; 0 otherwise\\
\end{tabular}

\begin{align}
\noindent \text{MSMP:}\notag \\
&z^*_L=\displaystyle\min_{\mathbf{w}} z(\mathbf{w}) \label{equation:zD}\\
&\indent\text{s.t.} \notag\\
& \indent  \displaystyle\sum_{i \in V} w_{i} \leq h \label{equation:budgetD1}\\
& \indent w_{i}\in \{0,1\} \indent &i& \in V \label{equation:binaryX}
\end{align}
where 
\begin{align}
&  z(\mathbf{w}) = \displaystyle\max_{\textbf{y},\textbf{u}} \displaystyle\sum_{r \in R}\sum_{i \in V}p_r u_{ir}  \label{equation:zA}\\
& \indent \text{\hspace{15pt}s.t.} \notag\\
& \indent\indent \displaystyle\sum_{i \in V} y_{i} \leq k \label{equation:budgetA1}\\
& \indent\indent y_{i}\leq 1-w_{i}  &i& \in V \label{equation:relateXY1}\\
& \indent\indent u_{ir} \leq  1- w_i  &i& \in V, r \in R \label{equation:seedSet1}\\
& \indent\indent u_{ir} \leq \sum_{j \in V} a_{ji}(\mathbf{w},r) y_j &i& \in V, r \in R \label{equation:active1}\\
& \indent\indent u_{ir}, y_{i} \in \{0,1\} &i& \in V, r \in R\label{equation:binary}
\end{align}

The objective \eqref{equation:zD} of the leader is to minimize the expected number of influenced nodes written explicitly in \eqref{equation:zA} as the follower's optimal objective value for the leader decision $\mathbf{w}$. Constraint \eqref{equation:budgetD1} and Constraint \eqref{equation:budgetA1} are the cardinality restrictions on the number of nodes protected and activated, respectively. Constraint \eqref{equation:relateXY1} and Constraint \eqref{equation:seedSet1} ensure that a protected node can neither be activated at the beginning nor be influenced in any scenario. For a node to be influenced in a scenario, it must be reachable from at least one seed node as stated in Constraint \eqref{equation:active1}. Finally, Constraint \eqref{equation:binaryX} and Constraint \eqref{equation:binary} are binary restrictions on the decision variables. Note that $R$ is defined as the set of all possible scenarios in \cite{kubra2019misinformation}, so that the distribution of the number of influenced nodes is identical to the one in the LT model. Since the size of $R$ is exhaustive for even very small networks, an approximation method is used to solve the follower's problem whereas the leader's problem is attacked with an heuristic. In this paper, we define $R$ as a live-arc scenario set of any size. Therefore, the spread distribution may be different than the one in the LT model.

The upper bound problem in the improved $x$-space algorithm, Algorithm \ref{algorithm:IXS}, corresponds to the follower's problem formulation in \eqref{equation:zA}--\eqref{equation:binary} and it needs to be updated and solved for each $\mathbf{w}^q$. To develop a lower bound problem formulation, we rearrange the model as follows. The term $a_{ji}(\mathbf{w},r)$ is a parameter which can be computed via a graph traversal algorithm for given $\mathbf{w}$, and it can be defined explicitly as
\begin{equation}
a_{ji}(\mathbf{w},r) = \left\{
        \begin{array}{ll}
            \prod_{k \in b_{ji}^r}(1-w_k) & \quad j \in a_i^r \\
            0 & \quad j \notin a_i^r,
        \end{array}
    \right.
\end{equation}
where $a_i^r$ denotes the set of nodes that can reach node $i$ via the arcs in scenario $r$ including $i$ itself. Parameter $b_{ji}^r$ represents the set of nodes on the path from $j$ to $i$ in scenario $r$ excluding $i$ and $j$, if such a path exists. In other words, there is an eligible path from $j$ to $i$ if none of the nodes between them is protected. By its definition, $a_{ji}(\mathbf{w},r)=1$ if $i=j$ or $j$ is the direct predecessor of $i$ in scenario $r$, which implies $b_{ji}^r=\varnothing$. Note that, there can be at most one path from one node to another in a scenario (a directed subgraph) due to the arc selection scheme of the live-arc representation for the LT model. Therefore, we are not concerned with a path selection decision. Let $\bar{a}_i^r$ denote the subset of $a_i^r$ excluding $i$ and its direct predecessor in scenario $r$, if it exists. To linearize the right-hand side of the \eqref{equation:active1}, we define a binary decision variable $m_{ji}^r$ for each $i \in V$, $j \in \bar{a}_i^r$ and $r\in R$, whose value is one if node $j$ is active and influences node $i$ in scenario $r$. In short $m_{ji}^r=y_j \prod_{k \in b_{ji}^r}(1-w_k)$ and constraints
\begin{align}
\hspace{2cm}& m_{ji}^r \leq y_j &i& \in V, j \in \bar{a}_i^r,r \in R \label{equation:m1}\\
&m_{ji}^r \leq  1-w_k &i& \in V, j \in \bar{a}_i^r, k \in b_{ji}^r,r \in R \label{equation:m2}\\
& m_{ji}^r \geq y_j -\sum_{k \in b_{ji}^r}w_k &i& \in V, j \in \bar{a}_i^r, r \in R  \label{equation:m3}
\end{align}
should be added to the lower level problem (LLP). Then, we replace Constraint \eqref{equation:active1} with
\begin{align}
&u_{ir} \leq \sum_{j \in a_i^r \setminus \bar{a}_i^r}y_j + \sum_{j \in \bar{a}_i^r} m_{ji}^r &i& \in V, r \in R. \label{equation:active_new}
\end{align}
Although constraint set \eqref{equation:m3} is required to define $m_{ji}^r$, it is easy to observe that, as a consequence of the optimization direction, i.e., maximization, the optimum objective value remains the same when this constraint set is relaxed. Hence, we relax it to obtain the final LLP; it consists of constraints \eqref{equation:zA}--\eqref{equation:seedSet1}, \eqref{equation:binary}, \eqref{equation:m1}--\eqref{equation:m2}, and \eqref{equation:active_new}. Note that it is also possible to formulate constraints \eqref{equation:active_new} as $u_{ir}\leq \sum_{j\in a_i^r} m_{ji}^r$, which would be still correct despite the fact that it leads to a larger MILP.

\subsection{Tailoring MSMP for the Improved $x$-space Algorithm}
\label{section:Implementation}
Recall that the lower bound (LB) problem \eqref{equation:LB problem} includes the convex combination constraints and the constraints that relate the lower and upper level variables. Let $\mathbf{x}=(\mathbf{y,u,m})$ denote the combined decision variables in the LLP of the MSMP. The LB problem formulation of the original $x$-space algorithm is provided below, where $ \mathcal{W}^\circ$ is defined by \eqref{equation:budgetD1} and \eqref{equation:binaryX}:

\begin{align}
l_q= \min_{\mathbf{w}\in \mathcal{W}^\circ}\, &\max_{\lambda,\mathbf{x}}  \sum_{r \in R}\sum_{i \in V}p_r u_{ir} \label{equation:MSMP-LB-obj}\\
& \hspace{10pt}\text{s.t.} \notag \\
& \hspace{10pt}\sum_{\tau\in J}\lambda_\tau=1 \\
& \hspace{10pt} \mathbf{x}-\sum_{\tau\in J}\lambda_\tau \mathbf{x}^\tau=0 \\
& \hspace{10pt} y_i \leq 1-w_i &i& \in V \label{equation:check-1}\\
& \hspace{10pt} u_{ir} \leq 1-w_i &i& \in V, r \in R \label{equation:check-2}\\
&  \hspace{10pt} m_{ji}^r \leq  1-w_k &i& \in V, j \in a_i^r, k \in b_{ji}^r,r \in R \label{equation:check-3}\\
& \hspace{10pt} \lambda \geq 0,\mathbf{x}\geq 0. \label{equation:MSMP-LB-nonnegative}
\end{align}

As can be noticed the formulation does not possess the special structure that we mention in Section \ref{subsection:new_LB}, i.e., $\mathbf{d}=\mathbf{1}$, $\mathbf{D}=\mathbf{I}$ and $n_1=n_2=n$. Despite this, the alternative LB problem formulations $LB_q'$ and $LB_q''$  can be derived for the MSMP as well. We first redefine the set $B(\mathbf{w},S_q)$, which allows us to obtain the particular form of $C_\tau$ necessary for these formulations.

\begin{proposition} For the MSMP,
$B(\mathbf{w},S_q)=\{\tau \in J : \exists i\in V, r\in R, w_i=1, u_{ir}^\tau=1 \}.$\label{Prop:MSMP blockers}
\end{proposition}

\begin{proof}
Recall that $B(\mathbf{w},S_q)$ is the index set of the follower solutions in $S_q$ that $\mathbf{w}$ blocks, i.e., renders infeasible. A follower solution becomes infeasible for $\mathbf{w}$ if one of the constraints in \eqref{equation:check-1}--\eqref{equation:check-3} is violated. Let $B_l(\mathbf{w},S_q)$ be the set of follower solutions that violate $l^{th}$ of these constraints. Then, $B(\mathbf{w},S_q)=\cup_{l=1}^3 B_l(\mathbf{w},S_q)$ by definition. We need to show that $B(\mathbf{w},S_q)=B_2(\mathbf{w},S_q)=\{\tau \in J : \exists i\in V, r\in R, w_i=1, u_{ir}^\tau=1 \}$. Consider $B_1(\mathbf{w},S_q)=\{\tau \in J : \exists i\in V, w_i=1, y_i^\tau=1 \}$. In an optimal follower solution, if $y_i=1$ then $u_{ir}=1$, for $r\in R$. This leads to the relation $B_1(\mathbf{w},S_q) \subset B_2(\mathbf{w},S_q)$, since the solutions in $S_q$ are optimal follower solutions except $\mathbf{x}^0$. Now consider $B_3(\mathbf{w},S_q)=\{\tau \in J \mid \exists i\in V, j \in a_i^r, k \in b_{ji}^r,r \in R, w_k=1, m_{ji}^\tau(r)=1 \}$. Choose $\tau \in B_3(\mathbf{w},S_q)$. If $m_{ji}^\tau(r)=1$, i.e., if node $j$ can influence node $i$ under scenario $r$ in solution $\tau$, then all the nodes on the path from $j$ to $i$ must be influenced as well. In other words, if $m_{ji}^\tau(r)=1$ then $u^\tau_{kr}=1$ for all $k \in b_{ji}^r$. Since also $w_k=1$ for some $k \in b_{ji}^r$, $\tau$ is contained in $B_2(\mathbf{w},S_q)$ implying that $B_3(\mathbf{w},S_q) \subset B_2(\mathbf{w},S_q)$.
\end{proof}

To summarize, any solution in conv($S_q$) that violates \eqref{equation:check-1} or \eqref{equation:check-3} also violates \eqref{equation:check-2}. This allows us to define the blocker set of a follower solution by using only constraint \eqref{equation:check-2}. Thus, we obtain $C_\tau=\{i\in V : \exists r\in R, u_{ir}^\tau=1\}$. Now, we can use $LB_q'$ or $LB_q''$ formulations developed in Section \ref{section:x-space}, since $C_\tau$ is explicitly defined. Note that for a fair comparison of the original and improved algorithms constraints \eqref{equation:check-1} and \eqref{equation:check-3} are removed from the lower bound problem formulation of the XS algorithm in \eqref{equation:MSMP-LB-obj}--\eqref{equation:MSMP-LB-nonnegative} as a result of Proposition \ref{Prop:MSMP blockers}.

\section{Computational Results} \label{section:Numerical Results}
In this section, we first examine the performances of the original $x$-space (XS) and improved $x$-space (IXS) algorithms, and compare them with a state-of-the-art method, MIX++ algorithm, developed in \cite{Fischettietal17} on the bilevel zero-one knapsack problem with interdiction (BKP) and bilevel maximum clique problem with interdiction (BCP). Then, we present the computational results related to the implementation of these algorithms on the MSMP described in Section \ref{section:MSMP}.

The XS and IXS algorithms are coded in C++ using Microsoft Visual Studio 2015. The experiments are carried out on a workstation with an Intel$\textregistered$ Xeon$\textregistered$ E5-2687W CPU, 3.10 GHz processor, and 64 GB RAM running within Microsoft Windows 7 Professional environment. The MIBLP solver, which is available at \url{https://msinnl.github.io/pages/bilevel.html}, is run on a virtual machine installed on the same workstation allocating 12 GB RAM with Ubuntu 20.04 operating system using CPLEX 12.7. The experiments involving XS and IXS algorithms are carried out with two MILP solvers: CPLEX 12.7 and Gurobi 8.0. We remark that in addition to the problems considered in this paper, the IXS algorithm can also be utilized for other well-known problem classes such as Shortest Path Interdiction Problem, Maximal Coverage Interdiction Problem, and $r$-Interdiction Median Problem.

\subsection{Results Obtained on the BKP and BCP Instances} \label{subsection:Comparison}

\cite{tang2016class} consider two types of interdiction problems, the bilevel knapsack problem (BKP), which is a min--max 0-1 knapsack problem, and the bilevel maximum clique problem (BCP), which is a min--max clique problem, to evaluate the performance of the improved $x$-space algorithm. In this section we compare the improved $x$-space algorithm with the original one as well as with MIX++  of \cite{Fischettietal17} by solving the instances given in \cite{tang2016class}. They can be reached at \url{http://jcsmith.people.clemson.edu/Test_Instances.html}.

The formulation of the BKP is given as
\begin{equation}
\min_{\mathbf{w}\in W} \max_{\mathbf{x}} \{p^T\mathbf{x} : a^T\mathbf{x} \leq b, \mathbf{x}\leq \mathbf{1}-\textbf{w} ,\mathbf{x} \in \{0,1\}^n\},
\label{equation:BKP}
\end{equation}
where
\begin{equation}
W=\{\mathbf{w}\in \{0,1\}^n : \mathbf{1}^T \mathbf{w} \leq k\}.
\end{equation}

Recall that it is necessary to define $C_\tau$, the blocker set of solution $\tau$, to develop the $LB_q''$ formulation of the improved algorithm. For the BKP, the set of solutions that $\mathbf{w}$ blocks is $B(\mathbf{w},S_q)=\{\tau \in J : \exists i \in  \{1,...,n\}, w_i=1, x_i^\tau=1\}$. Using this set, we can define $C_\tau$ as
\begin{equation}
C_\tau=\{i\in \{1,...,n\} : x_i^\tau=1\}.
\end{equation}

The second problem, BCP, is formulated as
\begin{align}
\min_{\mathbf{w}\in W} \, &\max \sum_{i \in V} x_i \label{equation:BCP1} \\
&\text{s.t. } x_i +x_j \leq 1 & (i,j)\in \bar{E} \\
&\hspace{20pt} x_i +x_j \leq 2-w_{ij}  & (i,j) \in E  \label{equation:BCP-interdiction}\\
&\hspace{20pt}x_i \in \{0,1\}\ &i\in V.\label{equation:BCP4}
\end{align}
Here $V$ and $E$ represent the set of vertices and edges, respectively, $\bar{E}=\{(i,j) : (i,j)\notin E \}$, and $W=\{ \mathbf{w} \in \{0,1\}^{|E|} : \sum_{(i,j) \in E} w_{ij} \leq k \}$.
Defining the set of blocked solutions according to the constraint set \eqref{equation:BCP-interdiction} as $B(\mathbf{w},S_q)=\{\tau \in J : \exists (i,j) \in E, w_{ij}=1, x^\tau_i+x^\tau_j=2\}$, gives rise to the blocker set
\begin{equation}
C_\tau=\{(i,j)\in E : x^\tau_i+x^\tau_j=2\}
\end{equation}
of edges as a result of the fact that the interdiction variables are related to the edges of the graph. In other words, if nodes $i$ and $j$ are included in the clique of solution $\mathbf{x}^\tau$, then interdicting the edge $(i,j)$ blocks that solution. Note that \cite{tang2016class} slightly modify the BCP formulation since the XS algorithm requires a special structure.

While obtaining the single-level LB formulation used in the $x$-space algorithm, it is required to dualize the inner problem in \eqref{equation:LB problem} and then linearize the resulting bilevel terms as explained in Section \ref{subsection:original_algorithm}. To this end, bounds on the dual variables of the LB problem are needed. For both BKP and BCP we compute the bounds as suggested in \cite{tang2016class} and use the same initial follower solution set $S_1$. It contains the optimal follower solutions for all possible unit vectors representing leader solutions, trivial follower solution $\mathbf{x}=\mathbf{0}$, and a heuristic feasible solution obtained via relaxing the integer variables of the follower's problem. We refer the reader to \cite{tang2016class} for details.

The results for 150 BKP benchmark instances with 15 different $(n,k)$ combinations are given in Table \ref{tab:BKP}. Columns XS-G and IXS-G are related to the results of our implementation of the $x$-space algorithm and improved $x$-space algorithm, respectively, which are obtained with Gurobi 8.0. Their counterparts with CPLEX 12.7 are given in columns XS-C and IXS-C. The last column titled MIX++ indicate the results given by the branch-and-cut solver of \cite{Fischettietal17} with its default settings where two types of intersection cuts are used along with a bilevel specific preprocessing procedure and follower upper bound cuts.
The results are presented in terms of the number of unsolved instances over 10 instances in the same $(n,k)$ setting, and the average CPU time for those instances with one-hour CPU time limit.

\begin{table}[h]
\fontsize{10pt}{11pt}\selectfont
  \centering
  \caption{Results obtained on BKP instances.}
    \begin{tabular}{cc|ccccc|ccccc}
\cline{3-12}          & \multicolumn{1}{r}{} & \multicolumn{5}{c|}{\# unsolved}      & \multicolumn{5}{c}{Time (sec.)} \\
    \hline
    \textit{n} & \textit{k} & XS-G & IXS-G & XS-C & IXS-C & MIX++ & XS-G & IXS-G & XS-C & IXS-C & \multicolumn{1}{c}{MIX++} \\
    \hline
    \multirow{3}[2]{*}{20} & 5     & 0     & 0     & 0     & 0     & 0     & 65.4  & 20.0  & 215.2 & 52.9  & 13.3 \\
          & 10    & 4     & 0     & 1     & 0     & 0     & 2486.4 & 213.4 & 1451.4 & 456.9 & 4.4 \\
          & 15    & 0     & 0     & 0     & 0     & 0     & 8.6   & 1.3   & 9.3   & 5.8   & 0.3 \\
    \hline
    \multicolumn{2}{c|}{Average} & 1.3   & 0.0   & 0.3   & 0.0   & 0.0   & 853.5 & 78.2  & 558.6 & 171.8 & 6.0 \\
    \hline
    \multirow{3}[2]{*}{22} & 6     & 0     & 0     & 4     & 0     & 0     & 1248.3 & 228.1 & 2724.9 & 610.9 & 21.5 \\
          & 11    & 8     & 1     & 8     & 1     & 0     & 3405.9 & 1376.7 & 3470.4 & 2074.3 & 6.7 \\
          & 17    & 0     & 0     & 0     & 0     & 0     & 20.9  & 1.8   & 12.4  & 7.3   & 0.3 \\
    \hline
    \multicolumn{2}{c|}{Average} & 2.7   & 0.3   & 4.0   & 0.3   & 0.0   & 1558.4 & 535.5 & 2069.3 & 897.5 & 9.5 \\
    \hline
    \multirow{3}[2]{*}{25} & 7     & 10    & 5     & 10    & 8     & 0     & 3601.9 & 2609.0 & 3601.6 & 3344.5 & 74.4 \\
          & 13    & 10    & 10    & 10    & 10    & 0     & 3601.3 & 3600.5 & 3600.6 & 3600.3 & 23.8 \\
          & 19    & 0     & 0     & 0     & 0     & 0     & 576.8 & 26.1  & 305.2 & 58.8  & 0.6 \\
    \hline
    \multicolumn{2}{c|}{Average} & 6.7   & 5.0   & 6.7   & 6.0   & 0.0   & 2593.3 & 2078.5 & 2502.5 & 2334.5 & 32.9 \\
    \hline
    \multirow{3}[2]{*}{28} & 7     & 10    & 10    & 10    & 10    & 0     & 3600.3 & 3600.6 & 3601.7 & 3601.2 & 191.7 \\
          & 14    & 10    & 10    & 10    & 10    & 0     & 3600.5 & 3600.4 & 3600.3 & 3600.5 & 63.3 \\
          & 21    & 8     & 0     & 5     & 0     & 0     & 3523.9 & 1042.2 & 3044.2 & 1194.4 & 1.2 \\
    \hline
    \multicolumn{2}{c|}{Average} & 9.3   & 6.7   & 8.3   & 6.7   & 0.0   & 3574.9 & 2747.7 & 3415.4 & 2798.7 & 85.4 \\
    \hline
    \multirow{3}[2]{*}{30} & 8     & 10    & 10    & 10    & 10    & 0     & 3600.3 & 3601.1 & 3600.6 & 3601.8 & 429.5 \\
          & 15    & 10    & 10    & 10    & 10    & 0     & 3600.6 & 3601.4 & 3600.3 & 3600.9 & 96.2 \\
          & 23    & 10    & 3     & 8     & 3     & 0     & 3600.6 & 1995.3 & 3472.0 & 1980.1 & 1.3 \\
    \hline
    \multicolumn{2}{c|}{Average} & 10.0  & 7.7   & 9.3   & 7.7   & 0.0   & 3600.5 & 3066.0 & 3557.6 & 3060.9 & 175.7 \\
    \hline
    \end{tabular}%
  \label{tab:BKP}%
\end{table}%

As can be observed, IXS performs better than XS in terms of the total number of unsolved instances with both MILP solvers. When the improved algorithm is used, the number of unsolved instances reduces from 90 to 59 with Gurobi and from 86 to 62 with CPLEX. Note that the number of unsolved instances for XS-C is less than the value reported in \cite{tang2016class}, which can be attributed to the MILP solver version (which is CPLEX 12.5 in that paper) or to hardware-related differences. When the CPU times are considered, the difference between the two algorithms becomes apparent especially for the $(n,k)$ settings in which most of the instances can be solved optimally within one hour. For example, the average CPU time for $(25,19)$ instances is reduced from 576.87 to 26.1 seconds with Gurobi, and from 305.2 to 58.8 seconds with CPLEX, whereas this number is equal to 1175 in \cite{tang2016class}. The number of unsolved instances corresponding to settings $(28,21)$ and $(30,23)$ decreases significantly for both solvers.

Nevertheless, MIX++ can solve all of the BKP instances optimally within one hour. This means that although the IXS algorithm is effective, it is outperformed by the state-of-the-art MIX++ method in solving BKP instances.

\begin{figure}[h!]
\centering
\begin{subfigure}{0.49\textwidth}
\includegraphics[width=1\linewidth]{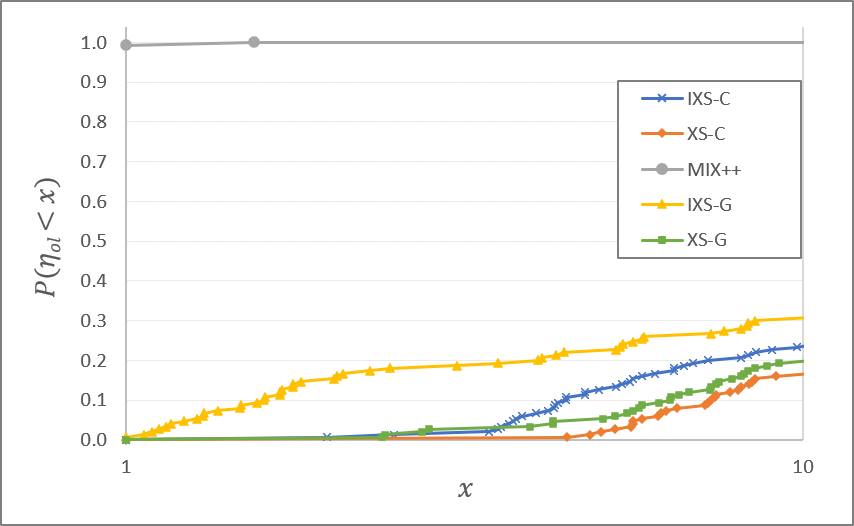}
\caption{$x\in[1,10]$}
\label{fig:BKP-1}
\end{subfigure}
\begin{subfigure}{0.49\textwidth}
\includegraphics[width=1\linewidth]{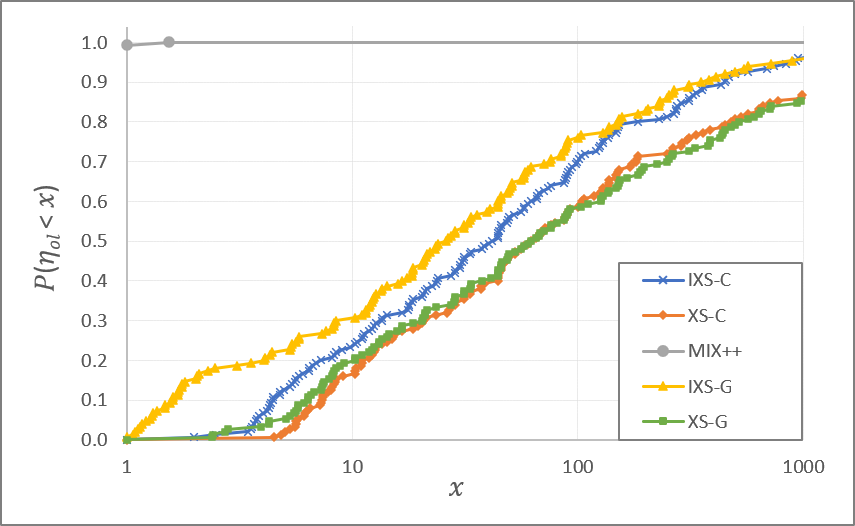}
\caption{$x\in[1,1000]$}
\label{fig:BKP-2}
\end{subfigure}
\caption{Comparison of the methods for BKP instances.}
\label{fig:BKP}
\end{figure}

Figure \ref{fig:BKP} shows the performance profiles of the five solution methods as described in \cite{dolan2002benchmarking} for benchmarking various optimization software. In this approach, a performance ratio $\eta_{o \ell}$ of method $\ell \in L$ on problem instance $o \in O$ is defined as the ratio of the solution time $\delta$ of $o$ with method $\ell$ to the minimum solution time for that problem instance, i.e., $\eta_{o \ell}=\delta_{o \ell}/\min_\ell \{\delta_{o\ell}\}$, when the performance measure of interest is the solution time $\delta$. The performance profile of each method $\ell \in L$ refers to the cumulative distribution function of $\eta_{o\ell}$. Following the approach in \cite{Fischettietal17}, we update the definition of $\eta_{o\ell}$ as follows to reduce the effect of the problem instances that can be solved in a very short amount of time (in seconds): $\eta_{o\ell}=(\delta_{o\ell} +1)/(\min_\ell \{\delta_{o\ell}\} +1)$. For those instances that cannot be solved by a given method within the time limit, $\delta_{o\ell}$ is set to 3600 seconds.

\begin{figure}[h!]
\centering
\includegraphics[width=0.5\linewidth]{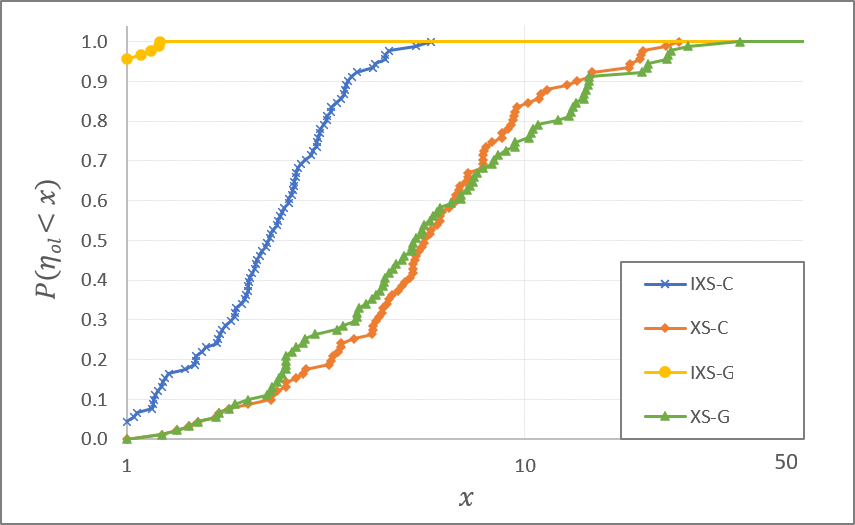}
\caption{Comparison of IXS and XS for BKP instances.}
\label{fig:BKPnew}
\end{figure}

Figures \ref{fig:BKP-1} and \ref{fig:BKP-2} use a logarithmic scale with different ranges, and it is clear that MIX++ is the best performer. Since it yields the minimum solution time in almost all BKP instances implying that $P(\eta_{o\ell}\leq 1)\approx 1$, it gives a straight line on the top of the chart. We can observe that IXS-G yields the next best setting, while XS-G and XS-C are very close to each other and yield the worst performance. The probability that IXS-G (IXS-C) solves a BKP problem at most 10 times slower than the best setting is around 30\% (23\%). Figure \ref{fig:BKPnew} shows the performance profiles of XS and IXS algorithms when MIX++ is not considered. IXS-G performs better than the others in $95\%$ of the instances while the remaining instances are solved faster by IXS-C. The probability that XS-G (XS-C) solution time is at most 10 times worse than the best method is $85\%$ ($75\%$).

The BCP instances in the benchmark data set are categorized into eight different $(n,d)$ combinations where $n$ denotes the number of nodes and $d$ denotes the network density. The number of interdicted edges is fixed to $k=\lceil {\frac{|E|}{4}} \rceil$ with $|E|$ representing the number of edges. There are 10 instances for each $(n,d)$ combination resulting in 80 instances in total. The results are provided in Table \ref{tab:BCP-penalty} in terms of the number of unsolved instances and CPU times. With the exception of the smallest instances with $(n,d)=(8,0.7)$, XS-C and XS-G fail to solve all the instances optimally within the time limit of one hour. The IXS algorithm can solve all the instances optimally within the time limit of one hour. As a matter of fact, it takes only a few seconds except the last setting with $(n,d)=(15,0.9)$, for which it requires less than two minutes on the average. MIX++ can also solve all of the instances optimally, although it yields a slightly larger average solution time than IXS-C and IXS-G. It is worth pointing out that MIX++ solution times that we obtain using its publicly available code are longer than the ones obtained via CPLEX 12.6.3 and reported in \cite{Fischettietal17}. Recall that we carry out the numerical experiments related to MIX++ using CPLEX 12.7 which is the newest version of CPLEX that the code works with at the time of our study.

\begin{table}[h!]
\fontsize{10pt}{11pt}\selectfont
  \centering
  \caption{Results obtained on BCP instances.}
    \begin{tabular}{cc|ccccc|ccccc}
\cline{3-12}          & \multicolumn{1}{r}{} & \multicolumn{5}{c|}{\# unsolved}      & \multicolumn{5}{c}{Time(sec.)} \\
    \hline
    $n$ & $d$ & XS-G & IXS-G & XS-C & IXS-C & MIX++ & XS-G & IXS-G & XS-C & IXS-C & MIX++ \\
    \hline
    \multirow{2}[2]{*}{8} & 0.7   & 0     & 0     & 0     & 0     & 0     & 1577.5 & 0.1   & 83.4  & 0.3   & 0.1 \\
          & 0.9   & 10    & 0     & 10    & 0     & 0     & 3600.4 & 0.1   & 3600.4 & 0.3   & 0.3 \\
    \hline
    \multicolumn{2}{c|}{Average} & 5.0   & 0.0   & 5.0   & 0.0   & 0.0   & 2588.9 & 0.1   & 1841.9 & 0.3   & 0.2 \\
    \hline
    \multirow{2}[2]{*}{10} & 0.7   & 10    & 0     & 10    & 0     & 0     & 3600.5 & 0.1   & 3600.5 & 0.4   & 0.3 \\
          & 0.9   & 10    & 0     & 10    & 0     & 0     & 3600.5 & 0.4   & 3600.4 & 0.7   & 1.3 \\
    \hline
    \multicolumn{2}{c|}{Average} & 10.0  & 0.0   & 10.0  & 0.0   & 0.0   & 3600.5 & 0.3   & 3600.5 & 0.5   & 0.8 \\
    \hline
    \multirow{2}[2]{*}{12} & 0.7   & 10    & 0     & 10    & 0     & 0     & 3600.6 & 0.3   & 3600.5 & 0.7   & 1.2 \\
          & 0.9   & 10    & 0     & 10    & 0     & 0     & 3600.4 & 2.2   & 3600.6 & 2.4   & 4.2 \\
    \hline
    \multicolumn{2}{c|}{Average} & 10.0  & 0.0   & 10.0  & 0.0   & 0.0   & 3600.5 & 1.3   & 3600.5 & 1.6   & 2.7 \\
    \hline
    \multirow{2}[2]{*}{15} & 0.7   & 10    & 0     & 10    & 0     & 0     & 3600.7 & 1.1   & 3600.6 & 1.8   & 5.4 \\
          & 0.9   & 10    & 0     & 10    & 0     & 0     & 3600.5 & 95.2  & 3600.6 & 113.1 & 168.5 \\
    \hline
    \multicolumn{2}{c|}{Average} & 10.0  & 0.0   & 10.0  & 0.0   & 0.0   & 3600.6 & 48.2  & 3600.6 & 57.5  & 87.0 \\
    \hline
    \end{tabular}%
\label{tab:BCP-penalty}
\end{table}%

\begin{figure}[h!]
\centering
\begin{subfigure}{0.49\textwidth}
\includegraphics[width=1\linewidth]{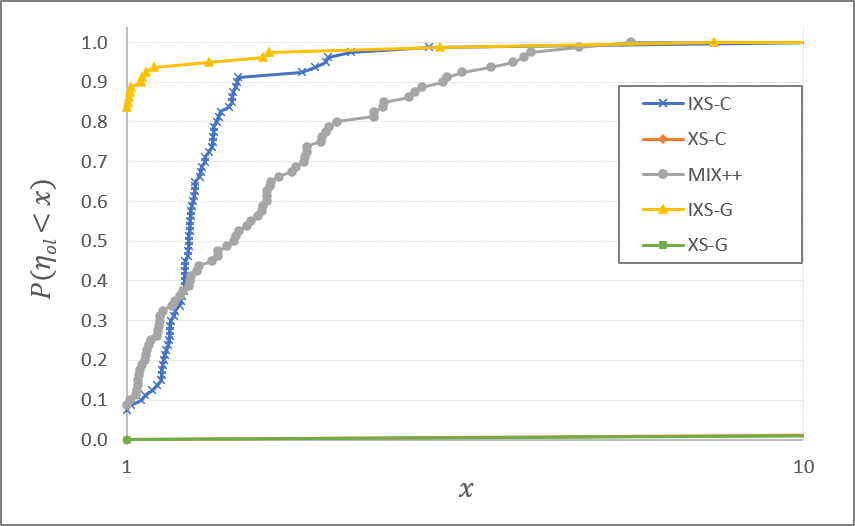}
\caption{$x\in[1,10]$}
\label{fig:BCP-1}
\end{subfigure}
\begin{subfigure}{0.49\textwidth}
\includegraphics[width=1\linewidth]{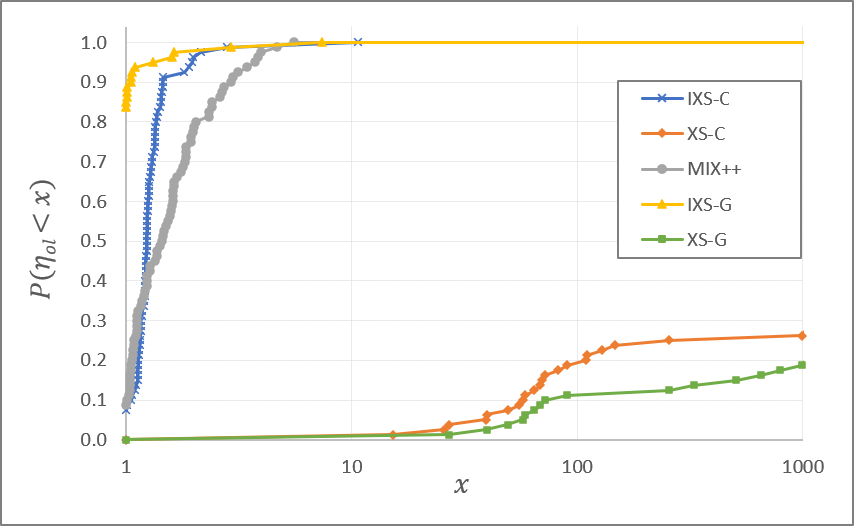}
\caption{$x\in[1,1000]$}
\label{fig:BCP-2}
\end{subfigure}
\caption{Comparison of the methods for BCP instances.}
\label{fig:BCP}
\end{figure}

The performance profiles are showed in Figure \ref{fig:BCP} on a logarithmic scale as before. They show that the proportion of instances for which IXS-G, IXS-C, and MIX++ have the minimum solution time are $80\%$, $10\%$, and $10\%$, respectively. The minimum performance ratios of XS-G and XS-C are 15 and 27, therefore their curves coincide with the horizontal axis in Figure \ref{fig:BCP-1}. IXS-C can solve $90\%$ of the instances within a CPU time of at most 1.5 times larger than the minimum time, whereas this ratio is 2.9 for MIX++. The proportion of instances for which the solution time is very close to the minimum (a performance ratio below 1.1) is larger for MIX++ compared to IXS-C. However, the proportion of instances with performance ratios greater than two is also larger for this method. In fact, IXS-C yields a ratio larger than two in only three of the instances and MIX++ does so in 17 of them.

Although the IXS algorithm performs better than XS for both interdiction problem types, the improvement it yields is much more apparent for the BCP. As a matter of fact, this improvement is due to the LB problems since the UB problems can be solved rather quickly. The average iteration times given in Table \ref{tab:Iter-times} show that the average time required to solve the LB problems is significantly smaller for IXS in both problem types. However, the number of iterations which determine the size of the final solution index set $J$ is still large with the BKP instances, as presented in Table \ref{tab:BKP-iter}. This implies that the success of IXS on BKP instances is due to the reduction in the average iteration times. For the BCP instances, both this reduction is larger and the number of iterations given in Table \ref{tab:BCP-iter} is significantly smaller. We interpret this situation as follows. The algorithm proceeds as long as a leader solution blocking all the follower solutions $\mathbf{x}^\tau$, $\tau \in J$, can be found. Finding such a solution is easier if the intersection of all $C_\tau$, $\tau \in J$ is larger. We observe that the average size of $C_\tau$ is significantly smaller for BCP where it is defined as the set of edges included in the maximum clique of the solution $\mathbf{x}^\tau$. This makes more difficult finding a leader solution blocking all $\mathbf{x}^\tau$ solutions, which causes the algorithm to quickly terminate. Regarding the reduction in average LB problem solution times, we attribute this result to the following difference between BKP and BCP. In the BKP, the interdiction relations are straightforward, i.e., in the form of $\mathbf{x} \leq \mathbf{1}-\mathbf{w}$ where the size of the upper and lower level variables are the same as mentioned in Section \ref{subsection:new_LB}, and all of the follower variables can be interdicted by the leader. However, in the BCP the leader interdicts some edges of the network and the follower selects nodes to form a maximum clique. The interdiction decisions indirectly cause some nodes to be excluded from the maximum clique. We can handle the situation by defining $C_\tau$ accordingly in the IXS algorithm without modifying the problem formulation.

When the performances of IXS and MIX++ are compared on the BKP and BCP instances, it is clear that IXS is outperformed on the former ones while it is the best performer in the latter ones. One possible reason is the size of the formulations used by the two methods. MIX++ solves the BKP instances using the bilevel formulation \eqref{equation:BKP} and the BCP instances using the formulation \eqref{equation:BCP1}--\eqref{equation:BCP4}. It is clear that for the same number of leader variables, the BCP formulation always includes more constraints and additionally these constraints have more nonzero coefficients as compared with BKP. On the other hand, IXS solves   whose number of variables and constraints depend only on the number of leader variables and the iteration number. In addition, IXS solves $LB_q''$ only if the greedy covering heuristic fails to find an optimal solution. This situation might give IXS an advantage on the BCP, in addition to smaller size of $C_\tau$ as explained above.

\subsection{Instance Generation and Results Obtained on the MSMP Instances} \label{subsection:MSMP Results}
The test instances for the MSMP are generated according to the Watts-Strogatz model which produces networks with small-world property, i.e., small distances between nodes and a relatively high clustering coefficient \citep{watts1998collective}. In the live-arc technique, scenarios (subgraphs) are sampled following a probability distribution determined via the arc weights. The weight parameters are generated uniformly in the $[0,1]$ interval and normalized to satisfy the restriction total incoming weight to a node cannot exceed one, as described in \cite{Kempe}. The live-arc scenario samples are generated using Latin Hypercube Sampling method since it is known that it yields smaller variance between samples, as validated in the preliminary experiments, and it is assumed that $p_r=1/|R|$, $r \in R$.

In all experiments of the IXS (XS) algorithm, the trivial follower solution selected in Step 0 of Algorithm \ref{algorithm:IXS} (Algorithm \ref{algorithm:XS}) is determined as $\mathbf{x}^0=\mathbf{0}$. Since any node in $V$ can be protected by the leader, this is the only solution that is feasible for all leader solutions. The initial set of the follower solutions is obtained by solving the follower's problem for randomly generated leader solutions. The preliminary experiments show that choosing leader solutions by simple heuristic methods instead of random selection does not have a significant impact on the results. For developing the UB problem formulations given in \eqref{equation:zA}--\eqref{equation:binary}, the value of $a_{ji}(\mathbf{w},r)$ is computed via a depth-first search on the reverse of the associated scenario subgraph, which is an in-tree, since at most one of the incoming arcs to a node can be live (active) in a scenario. For MIX++ the high-point relaxation problem is $\min_{\mathbf{w},\mathbf{x}} \big\{ \sum_{r \in R}\sum_{i \in V}u_{ir} :  \eqref{equation:budgetD1}, \eqref{equation:binaryX}, \eqref{equation:budgetA1}-\eqref{equation:seedSet1}, \eqref{equation:binary}, \eqref{equation:m1}--\eqref{equation:m2}, \eqref{equation:active_new} \big\} $.

Table \ref{tab:MSMP} displays the results over 200 test instances. There are 10 instances for each of the 20 $(n,R)$ combinations. The CPU time limit is one hour as before. Both XS and IXS algorithms (regardless of the MILP solver used) can solve all instances with 20 nodes within the time limit, with a significant difference in average solution times in favor of IXS. MIX++ fails to solve 15 instances with $n=20$, and leads to very long solution times. For $n\in \{25,30,35\}$, IXS-C and IXS-G can still solve all of the instances in one hour while XS-C fails to solve 68 and XS-G fails to solve 77 out of 120. It is observed that for these values of $n$, MIX++ yields smaller solution times than XS when $R=1$, while the total number of unsolved instances is 76. When the set of instances with $n=40$ is considered, XS-G and XS-C can solve only 9 and 10 of them, respectively, while IXS-G and IXS-C are able to solve 39 and 29 out of 40 instances optimally. When $n=40,R=1$, the average solution time of MIX++ is slightly smaller than IXS-C. However, it is outperformed by IXS-G, as is the case for all of the other problem sizes.

\begin{table}[h]
\fontsize{10pt}{11pt}\selectfont
  \centering
 \caption{Results obtained on MSMP instances.}
    \begin{tabular}{cc|ccccc|ccccc}
\cline{3-12}          & \multicolumn{1}{r}{} & \multicolumn{5}{c|}{\# unsolved}      & \multicolumn{5}{c}{CPU Time (sec.)} \\
    \hline
    $ n $     & $R$     & XS-G & IXS-G & XS-C & IXS-C & MIX++ & XS-G & IXS-G & XS-C & IXS-C & MIX++ \\
    \hline
    \multirow{4}[2]{*}{20} & 1    & 0 & 0     & 0     & 0     & 0     & 1.0   & 0.0   & 1.3   & 0.3   & 2.1 \\
          & 10    & 0     & 0     & 0     & 0     & 0     & 9.3   & 0.9   & 11.1  & 1.6   & 127.7 \\
          & 25    & 0     & 0     & 0     & 0     & 6     & 21.0  & 2.1   & 27.0  & 3.6   & 2309.5 \\
          & 50    & 0     & 0     & 0     & 0     & 9     & 46.3  & 7.6   & 52.3  & 8.7   & 3330.4 \\
    \hline
    \multicolumn{2}{c|}{Average} & 0.0   & 0.0   & 0.0   & 0.0   & 3.8   & 19.4  & 2.6   & 23.0  & 3.5   & 1442.4 \\
    \hline
    \multirow{4}[2]{*}{25} & 1     & 0     & 0     & 0     & 0     & 0     & 33.1  & 2.7   & 41.0  & 6.3   & 7.7 \\
          & 10    & 0     & 0     & 5     & 0     & 1     & 1090.6 & 19.2  & 1914.3 & 32.4  & 1060.6 \\
          & 25    & 3     & 0     & 5     & 0     & 6     & 2007.2 & 28.2  & 2525.3 & 54.0  & 2816.4 \\
          & 50    & 5     & 0     & 7     & 0     & 10    & 2288.9 & 40.0  & 2868.7 & 89.6  & 3600.1 \\
    \hline
    \multicolumn{2}{c|}{Average} & 2.0   & 0.0   & 4.3   & 0.0   & 4.3   & 1354.9 & 22.5  & 1837.3 & 45.6  & 1871.2 \\
    \hline
    \multirow{4}[2]{*}{30} & 1     & 0     & 0     & 0     & 0     & 0     & 786.1 & 29.7  & 365.3 & 50.2  & 53.9 \\
          & 10    & 10    & 0     & 10    & 0     & 10    & 3600.0 & 209.6 & 3604.4 & 436.1 & 3600.0 \\
          & 25    & 10    & 0     & 10    & 0     & 10    & 3600.0 & 237.4 & 3605.4 & 462.5 & 3600.2 \\
          & 50    & 10    & 0     & 10    & 0     & 10    & 3600.0 & 380.9 & 3604.7 & 706.0 & 3600.4 \\
    \hline
    \multicolumn{2}{c|}{Average} & 7.5   & 0.0   & 7.5   & 0.0   & 7.5   & 2896.5 & 214.4 & 2794.9 & 413.7 & 2713.6 \\
    \hline
    \multirow{4}[2]{*}{35} & 1     & 0     & 0     & 0     & 0     & 0     & 518.8 & 22.8  & 278.2 & 19.3  & 46.1 \\
          & 10    & 10    & 0     & 10    & 0     & 9     & 3600.0 & 314.3 & 3604.2 & 837.3 & 3534.0 \\
          & 25    & 10    & 0     & 10    & 0     & 10    & 3600.0 & 362.0 & 3604.6 & 878.3 & 3600.1 \\
          & 50    & 10    & 0     & 10    & 0     & 10    & 3600.0 & 418.8 & 3610.3 & 1099.1 & 3600.2 \\
    \hline
    \multicolumn{2}{c|}{Average} & 7.5   & 0.0   & 7.5   & 0.0   & 7.3   & 2829.7 & 279.5 & 2774.3 & 708.5 & 2695.1 \\
    \hline
    \multirow{4}[2]{*}{40} & 1     & 1     & 0     & 0     & 0     & 0     & 1813.4 & 67.7  & 1074.7 & 133.4 & 114.2 \\
          & 10    & 10    & 0     & 10    & 4     & 10    & 3600.0 & 1204.1 & 3602.7 & 2399.2 & 3600.1 \\
          & 25    & 10    & 0     & 10    & 4     & 10    & 3600.0 & 1175.7 & 3606.3 & 2192.3 & 3600.2 \\
          & 50    & 10    & 1     & 10    & 3     & 10    & 3600.0 & 1396.2 & 3607.2 & 2454.2 & 3600.3 \\
    \hline
    \multicolumn{2}{c|}{Average} & 7.8   & 0.3   & 7.5   & 2.8   & 7.5   & 3153.3 & 960.9 & 2972.7 & 1794.8 & 2728.7 \\
    \hline
    \end{tabular}%
  \label{tab:MSMP}%
\end{table}%

\begin{figure}[h!]
\centering
\begin{subfigure}{0.49\textwidth}
\includegraphics[width=1\linewidth]{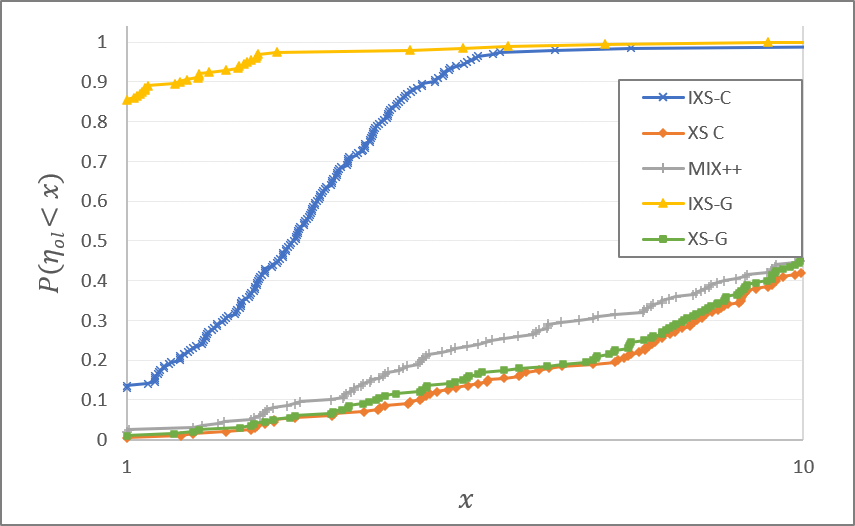}
\caption{$x\in[1,10]$}
\label{fig:MSMP-1}
\end{subfigure}
\begin{subfigure}{0.49\textwidth}
\includegraphics[width=1\linewidth]{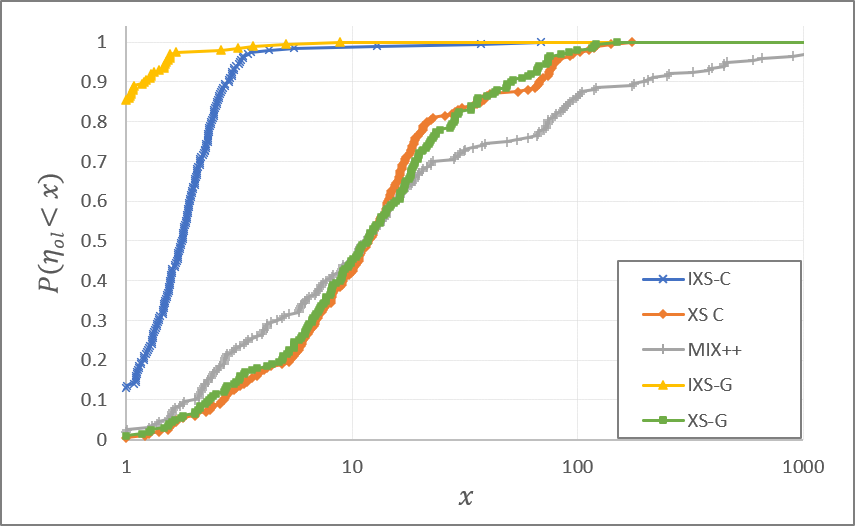}
\caption{$x\in[1,1000]$}
\label{fig:MSMP-2}
\end{subfigure}
\caption{Comparison of the methods for MSMP instances.}
\label{fig:MSMP}
\end{figure}

The performance profiles of each method over the MSMP instances are displayed in Figure \ref{fig:MSMP}. The intersections of the curves with the vertical axis indicate that IXS-G yields the smallest CPU time in $85\%$ of the instances and IXS-C needs the smallest time in $13\%$ of them. MIX++, XS-G and XS-C can only solve $2\%$, $1\%$ and $0.5\%$ of all instances respectively, faster than the other methods. The proportion of the instances that require at most five times more CPU time than the minimum solution time is $99.5\%$ for IXS-G, $98\%$ for IXS-C, $31\%$ for MIX++, $21\%$ for XS-G, and $19\%$ for XS-C. It can be seen in Figure \ref{fig:MSMP-2} that while MIX++ curve is above the ones of XS-C and XS-G for small values of $x$, it gets worse after $x=15$. In other words, while the probability of yielding a solution time that is at most 15 times worse than the minimum solution time is larger for MIX++ than XS-C and XS-G, the probability of yielding much worse performance ratios is also larger (e.g., $P(\eta_{o \ell}\leq 5)$ which is $31\%$ for MIX++ and $19\%$ for XS-C and $P(\eta_{o \ell}\leq 100)$ which is $86\%$ for MIX++, and $97\%$ for XS-C).

The results suggest that the IXS algorithm is less sensitive to the increase in the number of scenarios. Table \ref{tab:MSMP-iter} shows that the number of iterations usually increases as $R$ increases, which also increases the size of $C_\tau$, the set of nodes that are affected in any scenario in the solution $\mathbf{x}^\tau$ according to our observations. This result is in line with the findings of Section \ref{subsection:Comparison} which indicate that larger $C_\tau$ may lead to larger number of iterations, i.e., larger size of the final set $J$. We also observe that the size of $C_\tau$ is affected most when $R$ is increased from one to 10 and the average iteration times remain close for larger $R$. Although the number of iterations increases, IXS still performs well since the LB problems can be solved quickly.
The poor performance of MIX++ can be attributed to the size of the problems solved, since IXS makes it possible to discard a significant number of constraints as a result of Proposition \ref{Prop:MSMP blockers}. This result is also an indication for the necessity of algorithms that are tailored for the specific problem structure at hand.

\section{Conclusion} \label{section:Conclusion}
In this study, we improve an existing algorithm called $x$-space (XS) algorithm that is proposed for interdiction problems formulated as integer linear bilevel programming problems. The $x$-space algorithm iteratively solves lower and upper bound problems both of which are single level mixed-integer linear programs until the bounds converge to the same value. The solution time of the lower bound problem dominates the overall time required. Our methodology is based on developing an alternative lower bound problem formulation using the general features of the optimal solutions to this problem. The new formulation, which is similar to a maximum coverage problem, allows the use of a greedy coverage heuristic, which avoid solving a large number of lower bound mixed-integer linear programs exactly.

The improved $x$-space (IXS) algorithm has been firstly tested on various instances of two bilevel interdiction problems: bilevel knapsack problem with interdiction and bilevel maximum clique problem with interdiction. Then, it is implemented on a stochastic bilevel optimization problem on social networks (MSMP), by reformulating the problem and obtaining the problem specific blocker set definition. The IXS algorithm reduces the solution times significantly for all problem types considered as compared with the original XS algorithm. Furthermore, the computational results show that the IXS algorithm can handle the increase in the number of scenarios significantly better than the XS algorithm. A more general observation is that the improvement it yields becomes even more striking when the interdiction relation between the upper and lower level variables is not straightforward as is the case in the bilevel maximum clique problem with interdiction and misinformation spread minimization problem. When the performance of the IXS algorithm is compared with a state-of-the-art MIBLP solution algorithm MIX++, it is observed that although IXS is outperformed on the bilevel knapsack problem, it yields better solution times on the bilevel clique problem. Furthermore, it is clearly superior on the MSMP.

Searching for improvements in the lower bounding process via better formulations and faster algorithms remains as a challenge for new research endeavors. However, the quality of the upper bounds and the time spent for their computations are also important. Therefore, we believe that efforts in this direction might result in further enhancement regarding the overall efficiency of the algorithm.

\section*{Acknowledgements}
We thank the reviewers for their constructive comments that improved both the content and presentation of the paper. This study is partially supported by Bo\u{g}azi\c{c}i University Scientific Research Project under the Grant number: BAP 12745.

\section*{Appendix: Additional Computational Results}

\setcounter{table}{0}
\renewcommand{\thetable}{A\arabic{table}}

\begin{table}[htbp]
  \centering
  \fontsize{10pt}{11pt}\selectfont
  \caption{Average iteration times in seconds.}
    \begin{tabular}{ccccc}
    \hline
    Problem & XS-G & IXS-G & XS-C & IXS-C \\
    \hline
    BKP   & 0.31  & 0.08  & 0.41  & 0.15 \\
    BCP   & 0.50  & 0.02  & 0.48  & 0.03 \\
        MSMP  & 2.45  & 0.06  & 2.06  & 0.11 \\
    \hline
    \end{tabular} \label{tab:Iter-times}%
\end{table}%

\begin{table}[htbp]
\fontsize{10pt}{11pt}\selectfont
  \centering
  \caption{Number of iterations within the time limit for BKP instances.}
    \begin{tabular}{cc|cccc}
    \hline
    $ n $     & $R$     & XS-G & IXS-G & XS-C & IXS-C \\
    \hline
    \multirow{3}[2]{*}{20} & 5     & 1168.4 & 1242.7 & 1329.8 & 1133.7 \\
          & 10    & 4345.9 & 4276  & 3582.7 & 3355.5 \\
          & 15    & 217.4 & 315.4 & 291.8 & 259.0 \\
    \hline
    \multicolumn{2}{c|}{Average} & 1910.6 & 1944.7 & 1734.8 & 1582.7 \\
    \hline
    \multirow{3}[2]{*}{22} & 6     & 3808.8 & 4017.2 & 3575.7 & 3358.7 \\
          & 11    & 6957.3 & 10003.7 & 6870.4 & 7362.6 \\
          & 17    & 310.5 & 355.6 & 378.7 & 342.3 \\
    \hline
    \multicolumn{2}{c|}{Average} & 3692.2 & 4792.2 & 3608.3 & 3687.9 \\
    \hline
    \multirow{3}[2]{*}{25} & 7     & 10966.9 & 13088.7 & 4772.1 & 10576.5 \\
          & 13    & 11354.3 & 18553.7 & 9320.8 & 13563.2 \\
          & 19    & 1686.7 & 1727.9 & 1633.2 & 1566.0 \\
    \hline
    \multicolumn{2}{c|}{Average} & 8002.6 & 11123.4 & 5242.0 & 8568.6 \\
    \hline
    \multirow{3}[2]{*}{28} & 7     & 11326.4 & 21323.6 & 4621.6 & 14997.7 \\
          & 14    & 12431.1 & 37548.1 & 9762.9 & 21348.8 \\
          & 21    & 6602.8 & 9950.4 & 6405.1 & 6735.6 \\
    \hline
    \multicolumn{2}{c|}{Average} & 10120.1 & 22940.7 & 6929.9 & 14360.7 \\
    \hline
    \multirow{3}[2]{*}{30} & 8     & 11327.7 & 50446 & 6979.9 & 34652.6 \\
          & 15    & 11589.8 & 65263.7 & 10116.3 & 39036.6 \\
          & 23    & 7532.8 & 13981.0 & 8208.0  & 9222.6 \\
    \hline
    \multicolumn{2}{c|}{Average} & 10150.1 & 43230.2 & 8434.7 & 27637.3 \\
    \hline
    \end{tabular} \label{tab:BKP-iter}
\end{table}%

\begin{table}[htbp]
  \centering
  \fontsize{10pt}{11pt}\selectfont
  \caption{Number of iterations within the time limit for BCP instances.}
    \begin{tabular}{cc|cccc}
    \hline
     $ n $     & $R$     & XS-G & IXS-G & XS-C & IXS-C \\
    \hline
    \multirow{2}[2]{*}{8} & 0.7   & 5504.0  & 13.9  & 1239.7 & 13.7 \\
          & 0.9   & 10091.1 & 36.9  & 9557.4 & 35.7 \\
    \hline
    \multicolumn{2}{c|}{Average} & 7797.6 & 25.4  & 5398.6 & 24.7 \\
    \hline
    \multirow{2}[2]{*}{10} & 0.7   & 9065.6 & 29.9  & 8543.1 & 27.2 \\
          & 0.9   & 7729.2 & 92.2  & 7683.4 & 93.9 \\
    \hline
    \multicolumn{2}{c|}{Average} & 8397.4 & 61.1  & 8113.3 & 60.6 \\
    \hline
    \multirow{2}[2]{*}{12} & 0.7   & 7198.1 & 67.6  & 6768.7 & 64.3 \\
          & 0.9   & 6105.4 & 241.2 & 6365.4 & 237.2 \\
    \hline
    \multicolumn{2}{c|}{Average} & 6651.8 & 154.4 & 6567.1 & 150.8 \\
    \hline
    \multirow{2}[2]{*}{15} & 0.7   & 5603.0  & 175.8 & 5432.2 & 165.3 \\
          & 0.9   & 4272.5 & 791.0   & 4768.2 & 767.4 \\
    \hline
    \multicolumn{2}{c|}{Average} & 4937.8 & 483.4 & 5100.2 & 466.4 \\
    \hline
    \end{tabular} \label{tab:BCP-iter}
\end{table}%

\begin{table}[htbp]
  \centering
  \fontsize{10pt}{11pt}\selectfont
  \caption{Number of iterations within the time limit for MSMP instances.}
    \begin{tabular}{cc|cccc}
    \hline
        $ n $     & $R$     & XS-G & IXS-G & XS-C & IXS-C \\
    \hline
    \multirow{4}[2]{*}{20} & 1     & 41.7  & 35.0    & 50.1  & 34.1 \\
          & 10    & 107.3 & 100.3 & 105.7 & 100.4 \\
          & 25    & 126.8 & 116.8 & 124.2 & 117.6 \\
          & 50    & 136.4 & 127.3 & 134.3 & 128.4 \\
    \hline
    \multicolumn{2}{c|}{Average} & 103.1 & 94.9  & 103.6 & 95.1 \\
    \hline
    \multirow{4}[2]{*}{25} & 1     & 260.8 & 218.6 & 354.7 & 184.2 \\
          & 10    & 1174.5 & 1297.6 & 1222.1 & 1318.0 \\
          & 25    & 1340.5 & 1467.4 & 1055.6 & 1482.8 \\
          & 50    & 1397.4 & 1524.4 & 920.7 & 1547.7 \\
    \hline
    \multicolumn{2}{c|}{Average} & 1043.3 & 1127.0 & 888.3 & 1133.2 \\
    \hline
    \multirow{4}[2]{*}{30} & 1     & 1119.4 & 747.9 & 1115.2 & 665.3 \\
          & 10    & 3366.8 & 3401.7 & 1913.1 & 3415.8 \\
          & 25    & 2375.1 & 3734.2 & 1108.9 & 3747.1 \\
          & 50    & 1829.4 & 3861.3 & 829.8 & 3871.2 \\
    \hline
    \multicolumn{2}{c|}{Average} & 2172.7 & 2936.3 & 1241.8 & 2924.9 \\
    \hline
    \multirow{4}[2]{*}{35} & 1     & 759.3 & 611.8 & 821.7 & 421.3 \\
          & 10    & 3300.2 & 5070.2 & 1907.5 & 5135.5 \\
          & 25    & 2458.3 & 5628.7 & 992.4 & 5684.9 \\
          & 50    & 1879.8 & 5796.0  & 747.5 & 5845.1 \\
    \hline
    \multicolumn{2}{c|}{Average} & 2099.4 & 4276.7 & 1117.3 & 4271.7 \\
    \hline
    \multirow{4}[2]{*}{40} & 1     & 1530.5 & 1078.4 & 1444.6 & 757.5 \\
          & 10    & 3257.0  & 8721.3 & 1758.3 & 8760.0 \\
          & 25    & 2181.9 & 9367.4 & 924.8 & 9387.0 \\
          & 50    & 1671.9 & 9485.3 & 633.0   & 9465.6 \\
    \hline
    \multicolumn{2}{c|}{Average} & 2160.3 & 7163.1 & 1190.2 & 7092.5 \\
    \hline
    \end{tabular} \label{tab:MSMP-iter}
\end{table}

\bibliographystyle{apalike}


\end{document}